\documentclass{amsart}
\usepackage{amsfonts,amssymb,amsmath,a4wide}
\usepackage[british]{babel}
\newcommand{\qedbox}{ \fbox{}}

\newtheorem{teo}{Theorem}[section]
\newtheorem{lema}{Lemma}[section]
\newtheorem{pro}{Proposition}[section]
\newtheorem{defi}{Definition}[section]\newtheorem{rema}{Remark}[section]
\newtheorem{coro}{Corollary}[section]
\numberwithin{equation}{section}

\def \V{\nu}

\def \z{\zeta}

\def \l4{[[\Lambda^{0,4}]]}
\def \ll3{[[\Lambda^{0,3}]]}

\def \g{\mathfrak{g}}
\def \SS{\slash \hspace{-2.3mm}S}                                                                                               
\newcommand{\Real}{\mathbb{R}}                                                                                                  
\newcommand{\Cl}{{Cl}}                                                                                                      

\newcommand{\solie}{\mathfrak{so}}

\let\phi=\varphi
\let\star=*
\newcommand{\vph}{\varphi}                                                                                                      
\newcommand{\h}{\frac{1}{2}}                                                                                                    
\newcommand{\sh}{{\textstyle \frac{1}{2}}}                                                                                      
\def\hook{\mathbin{\hbox to 6pt{\vrule height0.4pt width5pt depth0pt \kern-.4pt \vrule height6pt width0.4pt depth0pt\hss}}}     
\let\lrcorner=\hook

\def\sideremark#1{\ifvmode\leavevmode\fi\vadjust{\vbox to0pt{\vss
 \hbox to 0pt{\hskip\hsize\hskip1em
 \vbox{\hsize2.5cm\tiny\raggedright\pretolerance10000
 \noindent #1\hfill}\hss}\vbox to8pt{\vfil}\vss}}}%

\begin{document}\larger[2]
\title[]{On algebraic torsion forms and their spin holonomy algebras}
\author[N.Bernhardt]{Niels Bernhardt}
\address[N.Bernhardt]{Department of Mathematics, University of Auckland, Private Bag 92019, Auckland, New Zealand}
\email{niels.bernhardt@math.auckland.ac.nz}
\author[P.-A. Nagy]{Paul-Andi Nagy}
\address[P.-A. Nagy]{Department of Mathematics, University of Auckland, Private Bag 92019, Auckland, New Zealand}
\email{nagy@math.auckland.ac.nz}
\date{\today}

\begin{abstract}
We study holonomy algebras generated by an algebraic element of the Clifford algebra, or equivalently, the holonomy algebras of certain spin connections in flat space. We provide
series of examples in arbitrary dimensions and establish general properties of the holonomy algebras under some mild conditions on the generating element. We show that the first
non-standard situation to look at appears in dimension $8$ and concerns self-dual $4$-forms. In this case complete structure results are obtained.

\end{abstract}

\maketitle
\tableofcontents

\section{Introduction}
Let $(M^n,g)$ be a Riemannian spin manifold, with spinor bundle to be denoted by $\SS$. For any differential form $T$ on $M$,
not necessarily of pure degree, one can form the linear connection $\nabla^T$ on $\SS$ by setting
\begin{equation*}
\nabla_X^T \psi=\nabla_X \psi+(X \lrcorner T) \psi
\end{equation*}
whenever $\psi$ belongs to $\Gamma(\SS)$ and $X$ is in $TM$. Here $\nabla$ is the connection induced by the Levi-Civita connection on the spinor bundle $\SS$. This can be thought of as the spin
analogue of a connection with non-trivial torsion on the tangent bundle of $M$. A special case is when $T$ is actually
a $3$-form in which case $\nabla^T$ is the induced spin connection of a connection with torsion on $TM$. In low dimensions, ranging from $6$ to $8$, parallel
spinors w.r.t to a connection with $3$-form torsion are nowadays rather well understood in terms of geometric structures on the tangent bundle to the manifold \cite{brya,ivan}
and extensive effort toward their classification has been made \cite{mafe,GrayH,FI}.
This has been also studied in connection with the so-called Strominger's type II string equations \cite{stro}, \cite{AFNP}. Another special case, which no longer reflects the presence
of a particular connection at the level of the tangent bundle of $M$, is when $T$ consists of forms of degree $3$ and $4$, the latter being termed fluxes in physics literature (see \cite{papa}
and references therein). In all the above mentioned cases one of the issues to understand is under which conditions $\nabla^T$ admits parallel spinors, therefore one looks,
more generally, at the holonomy representation of $\nabla^T$. This is because of the well-known fact \cite{Besse} that the existence of a parallel spinor field is equivalent with the spinor being fixed
by the holonomy representation at a point. \\
In this paper we shall study the holonomy of the connection $\nabla^T$ in the flat case, when moreover $T$ is assumed to have constant coefficients. This is the simplest geometric case one could think
of but already raises some interesting and quite difficult algebraic questions. We set
\begin{defi} \label{maindf}Let $(V^n,  \langle \cdot, \cdot \rangle)$  be an Euclidean vector space and let $Cl_n(V)$ be its Clifford
algebra. Then:
\begin{itemize}
\item[(i)] the fix algebra of $T$ in $Cl_n(V)$ is the Lie-sub-algebra $\mathfrak{g}_T^{\star}$ of $Cl_n(V)$ generated by $\{X \lrcorner T: X \in V \}$.
\item[(ii)] the holonomy algebra of some $T$ in $Cl_n(V)$ is given as $\mathfrak{h}_T^{\star}=[\mathfrak{g}_T^{\star},\mathfrak{g}_T^{\star}]$.
\end{itemize}
\end{defi}
This is motivated by the observation \cite{holo} that in the flat case the holonomy algebra of the spin connection $\nabla^T$ equals $\mathfrak{h}_T^{\star}$. When $\mathfrak{g}_T^{\star}$ is {\it{perfect}},
that is $\mathfrak{g}_T^{\star}=[\mathfrak{g}_T^{\star},\mathfrak{g}_T^{\star}]$, the two algebras above coincide and in this respect the fix algebra
$\mathfrak{g}_T^{\star}$ appears to be a very useful intermediary object for establishing structure results, although it seems to lack of further geometric content. For the special $3$-forms
complete structure results concerning the fix and holonomy algebras have been obtained in \cite{holo}.\\
Our paper is organised as follows. In section 2 we review a number of elementary facts concerning Clifford algebras and their representations, with accent put on the different
phenomena appearing in some arithmetic series of dimensions. In section 3 we start our study of holonomy algebras by determining - under some mild assumptions on the generating
element - the model algebra those are contained in. We also establish a number of useful general properties, like semisimplicity.
Further on, we investigate the space of the so-called fixed spinors which for some $T$ in $Cl_n(V)$ is defined as
\begin{equation*} \label{spinfixdef}
Z_T=\{\psi \in \SS: (X \lrcorner T)\psi=0 \ \mbox{for all} \ X \ \mbox{in} \ V\}
\end{equation*}
where $\SS$ is an irreducible $Cl_n(V)$ module. We provide first order information about these spaces and also
discuss some simple examples. The section ends with giving a necessary condition for certain holonomy algebras to be perfect, namely
\begin{teo} \label{mainperf}Let $T$ in $Cl_n^0 \cap Cl_n^{+}$, where $n \equiv 0$ (mod 4) satisfy $T^t=T$. If $Z_T=(0)$ then $\mathfrak{g}_T^{\star}$ is a perfect Lie algebra.
\end{teo}
Section 4 describes situations where the holonomy algebras can be directly computed and provides series of useful, in hindsight, examples. Elements of the Clifford algebra which are
being looked at are unipotent and squares of spinors, which actually give idempotents. In the latter situation, the dimension (mod 8) of the underlying vector space appears to lead to
very different results. More precisely
\begin{teo} \label{main sq}Let $V$ be an Euclidean vector space with volume form $\nu$ and let $T$ belong to $Cl_n(V)$. Then:
\begin{itemize}
\item[(i)] if $n \equiv 0$ (mod 4) and $T$ in $Cl_n^0 \cap Cl_n^{+}$ is unipotent, that is $T^2=1+\nu$, then its holonomy algebra is isomorphic to $\mathfrak{so}(n,1)$ and $\mathfrak{g}_T^{\star}$
is perfect.
\item[(ii)] if $T$ is the square of a spinor when $n \equiv 7$ (mod 8), then $\mathfrak{g}_T^{\star}$ is abelian, in particular the holonomy algebra vanishes.
\item[(iii)] if $T$ is the square of a positive spinor when $n \equiv 0$ (mod 8), then the fix algebra of $T$ is perfect and its holonomy algebra is isomorphic to $\mathfrak{so}(n,1)$.
\end{itemize}
\end{teo}
This essentially exploits specific features of the powerful squaring construction for spinors \cite{laws}, \cite{dado}. Note that the perfecteness of the fix algebra in (i) of Theorem \ref{main sq} follows
directly from Theorem \ref{mainperf} wheareas in the case of (iii) it does not, since the set of fixed spinors is no longer trivial.
Section 5 forms the core of the present paper and gives a complete classification of holonomy algebras generated by self-dual $4$-forms in dimension $8$. More precisely, we show
\begin{teo} \label{main8} Let $V$ be an oriented Euclidean vector space of dimension $8$ and let $T$ be a self-dual four form on $V$. The fix algebra of $T$ is perfect and
its holonomy algebra is isomorphic to the Lie algebra $\mathfrak{so}(8,8-dim_{\mathbb{R}}Z_T)$ exception made of the cases when
\begin{itemize}
\item[(i)] $T$ is a unipotent element
\item[(ii)] $dim_{\mathbb{R}}Z_T=6$
\end{itemize}
where the holonomy algebras are $\mathfrak{so}(8,1)$ and $\mathfrak{so}(6,2)$ respectively.
\end{teo}
The proof, the case when  $dim_{\mathbb{R}}Z_T=6$ excepted, uses the splitting of the space of two forms which can obtained once given a self-dual four form.
This is combined with the observation that raising the generating element to any odd Clifford power leaves the initial fix algebra unchanged. In the last section of
the paper, we treat directly the special case appearing in (ii) of Theorem \ref{main8} using the one to one correspondence \cite{HL} between the existence of a such a form
and that of an $SU(4)$-structure on our vector space. The paper ends by an appendix, containing the elementary though lengthy proof of a technical Lemma.

\section{Preliminaries}
This section is mainly intented to recall a number of facts concerning Clifford algebras and spinors, which we shall constantly use in what follows. A thorough account of all theses notions can
be found in \cite{laws}.
\subsection{Clifford algebras}

Let $V$ be an $n$-dimensional vector space over $\mathbb{R}$ equipped with a scalar product, to be denoted by $\langle \cdot, \cdot \rangle$. We shall denote by $\Cl_n(V)$ the Clifford
algebra associated with $(V, \langle \cdot, \cdot \rangle)$, and if there is no ambiguity on the vector space used we shall simply write $\Cl_n$ for $\Cl_n(V)$. We recall that $\Cl_n(V)$ can be
given the structure of an algebra, with multiplication denoted by $" \cdot ": \Cl_n(V) \to \Cl_n(V)$ and satisfies
\begin{equation} \label{clifford1form}
    e \cdot \vph = e \wedge \vph - e \hook \vph, \qquad \vph \cdot e = (-1)^k (e \wedge \vph + e \hook \vph)
\end{equation}
whenever $e$ belongs to $\Lambda^1(V)$ and $\varphi$ is in $\Lambda^{k}(V)$, although this notation will no longer be used in what follows. Here and henceforth we will identify $1$-forms and vectors
via the given scalar product $\langle \cdot, \cdot \rangle$. There is a canonical isomorphism between the space $\Lambda^*(V)$ of forms in $V$ and the Clifford algebra $\Cl_n(V)$ having the
following property. Let $L: \Cl_n(V) \to \Cl_n(V)$ be defined by
\begin{equation*} \label{loperator}
    L(\vph) = \sum \limits_{i=1}^{n} e_i \, \vph \, e_i,
\end{equation*}
whenever $\vph \in \Cl_n(V)$ and for some orthonormal basis $\{e_i\}$, $1 \leq i \leq p$ in $V$. Then the eigenspaces of $L$ are the canonical images of $\Lambda^k(V)$:
\begin{equation} \label{leigenspaces}
    L|_{\Lambda^k}=(-1)^k (2k-n) 1_{\Lambda^k}.
\end{equation}
Any Clifford algebra comes with two involutions, the first being the transposition map
$(\, )^t: \Cl_n(V) \to \Cl_n(V)$
defined by
\begin{equation*} \label{cliffordtranspose}
    (e_1 \, e_2 \, ... \, e_{k-1} \, e_k)^t = e_k \, e_{k-1} \, ... \, e_2 \, e_1,
\end{equation*}
for some orthonormal frame $\{e_i,1 \le i\le n \}$. Note however that the transpose is frame independent and therefore extends to an anti-automorphism of $Cl_n(V)$, i.e.
\begin{equation} \label{transposerule}
    (\vph_1 \, \vph_2 )^t = \vph_2^t \, \vph_1^t
\end{equation}
for all $\vph_1, \vph_2$ in $\Cl_n(V)$. The second involution $\alpha: Cl_n(V) \to Cl_n(V)$ results from extending $-1_{V}$ to an automorphism of the algebra $Cl_n(V)$, in the sense that
\begin{equation*}
    \alpha(\vph_1 \, \vph_2) = \alpha(\vph_1) \, \alpha(\vph_2),
\end{equation*}
where $\vph_1, \vph_2$ are in $\Cl_n(V)$. Since $\alpha$ is an involution it can also be used to obtain a splitting
\begin{equation*}
    \Cl_n(V) = \Cl_n^{0}(V) \oplus \Cl_n^{1}(V)
\end{equation*}
into the $\pm$-eigenspaces of $\alpha$. Note this corresponds to the splitting of $\Lambda^*(V)$ into even and respectively odd degree forms. The vector space $\Cl_n(V)$ inherits from $V$ a scalar
product, still to be denoted by $\langle \cdot, \cdot \rangle$ and having the property that
\begin{equation} \label{scprod}
\begin{split}
    & \langle q  \vph_1, \phi_2 \rangle = \langle \vph_1, \alpha(q^{t}) \vph_2 \rangle \\
    & \langle \vph_1  q, \phi_2 \rangle = \langle \vph_1, \vph_2  \alpha(q^{t}) \rangle
\end{split}
\end{equation}
whenever $q, \vph_1, \vph_2$ belong to $\Cl_n(V)$. Let us assume now that $V$ is oriented by $\nu$ in $\Lambda^n(V)$ such that for an oriented frame $\{e_k, 1 \le k \le n \}$ this is given as
$\nu = e_1 \, \ldots \, e_n$. Then it is easy to check that
\begin{equation} \label{vsquare}
    \nu^2=(-1)^{\frac{n(n+1)}{2}}, \qquad  \nu^t=(-1)^{\frac{n(n-1)}{2}}\nu.
\end{equation}
Now the Hodge star operator  $*: \Lambda^k(V) \to \Lambda^{n-k}(V)$ is defined by $\alpha \wedge *\beta = \langle \alpha, \beta \rangle \nu$, for all $\alpha, \beta$ in $\Lambda^*(V)$
and relates to Clifford multiplication with $\nu$ by
\begin{equation} \label{cliffordhodge}
    * \vph = (-1)^{\frac{k}{2}(k+1)} \vph  \nu = (-1)^{\frac{k}{2}(2n-k+1)} \nu  \vph,
\end{equation}
for $\vph$ in $\Lambda^k(V) \subset \Cl_n(V)$. Moreover, we have
\begin{equation} \label{alphaeven}
\begin{array}{ll}
    \vph \, \nu = \nu \, \vph,            & \mbox{for all } \ \phi \ \mbox{in}\  \Cl_n^0, \cr
    \vph \, \nu = (-1)^{n+1} \nu \, \phi, & \mbox{for all } \ \phi \ \mbox{in} \  \Cl_n^1.
\end{array}
\end{equation}
In particular, when $n \equiv 1$ (mod 2), the volume element $\nu$ belongs to the center of the Clifford algebra $\Cl_n$.

If $n \equiv 0$ (mod 4) then $\nu^2 = 1$ whence the Hodge star operator, realised as in \eqref{cliffordhodge}, provides a decomposition of the Clifford algebra into self-dual and anti-self-dual elements:
\begin{equation} \label{hodgeeven}
    \Cl_n(V) = \Cl^{+}_n(V) \oplus \Cl^{-}_n(V),
\end{equation}
where $\V \vph =\pm  \vph $ whenever $\vph$ belongs to $\Cl^{\pm}_n(V)$.

\subsection{The space of spinors}
We need to recall some elementary facts about spinors.  Let $\SS$ be an irreducible $\Cl_n(V)$-module. We shall call $\SS$ the space of spinors and elements $\psi \in \SS$ spinors and
denote by $\mu: \Cl_n(V) \to {\rm End}(\SS)$ the Clifford multiplication acting on $\SS$. On $\SS$ we have the usual scalar product $\langle \psi_1, \psi_2 \rangle$ for two
spinors $\psi_1, \psi_2 \in \SS$ 
which has the following property
\begin{equation} \label{spinorscalarprod}
    \langle \vph  \psi_1, \psi_2 \rangle = \langle \psi_1, \alpha(\vph^t) \psi_2 \rangle,
\end{equation}
for all $\vph$ in $\Cl_n$ and $\psi_1, \psi_2$ in $\SS$. Recalling that
\begin{equation}
\alpha(\vph^t) = (-1)^{\frac{k}{2}(k+1)} \, \vph, \qquad \vph \in \Lambda^k(V) \subset \Cl_n(V).
\end{equation}
it follows that the Clifford multiplication operator $\mu_{\vph} :\SS \to \SS$, with an element $\vph$ of $\Lambda^k(V)$, is symmetric when $k \equiv 0,3$ (mod 4) and anti-symmetric
when $k \equiv 1,2$ (mod 4). Now when $n \equiv 0, 3$ (mod 4) the volume form $\nu$ squares to $1$. If $n \equiv 0$ (mod 4) this allows splitting the irreducible, real Clifford module $\SS$
as $\SS=\SS^{+} \oplus \SS^{-}$, where $\nu$ acts as $\pm 1$ on $\SS^{\pm}$. Peculiar to the case when $n \equiv 3$ (mod 4) is the fact that any irreducible, real, Clifford module
$\SS$ has either $\nu \psi = -\psi$ for all $\psi$ in $\SS$ or $\nu \psi=\psi$ for all $\psi$ in $\SS$. Both possibilities can occur and produce different $\Cl_n$-representations. As a
convention, in what follows we shall always work with the latter representation. Let us also mention that when $n \equiv 3$ (mod 4) we have that $\nu \phi = \phi \nu$ for all $\phi$ in
$\Cl_n$ and that $\alpha$ interchanges $\Cl_n^{+}$ and $\Cl_n^{-}$, that is realises an isomorphism $\alpha :\Cl_n^{+} \to \Cl_n^{-}$. For it will be used constantly in what follows we
also recall the following stability Lemma.
\begin{lema} \label{stability}
The following stability conditions hold:
    \begin{itemize}
    \item[(i)]  when $n \equiv 0$ (mod 4),
                \begin{equation*}
                    \begin{array}{ll}
                        \vph  \SS^+ \subseteq \SS^+, & \forall \vph \in \Cl^0_n(V) \cap \Cl^+_n(V) \cr
                        \vph  \SS^+ = 0,             & \forall \vph \in \Cl^0_n(V) \cap \Cl^-_n(V) \cr
                        \vph  \SS^+ = 0,             & \forall \vph \in \Cl^1_n(V) \cap \Cl^+_n(V) \cr
                        \vph  \SS^+ \subseteq \SS^+, & \forall \vph \in \Cl^1_n(V) \cap \Cl^-_n(V),
                    \end{array}
                \end{equation*}
    \item[(ii)] while for $n \equiv 3 $ (mod 4) we have $\Cl_n^{+}\SS \subseteq \SS, \Cl_n^{-} \SS=0$.
\end{itemize}
\end{lema}

The proof, which is left to the reader, follows from the above properties of Clifford multiplication with $\nu$. Similar statements can be easily made on $\SS^{-}$ when $n \equiv 0$
(mod 4). We end this section by recalling two more well known facts, with proofs given for the sake of completeness.

\begin{lema} \label{prep}
If $\z$ in $\Cl_n$ satisfies $[\z, \Lambda^2(V)] = 0$ then $\z$ belongs to $(1, \nu)$.
\end{lema}
\begin{proof}
Since $\z XY = XY \z$ for all $X,Y$ in $V$ it follows that $Y(X\z X)Y=\vert X \vert^2 \vert Y \vert^2 \z$ and further, after a double tracing $L^2 \z=n^2 \z$. The eigenvalues of $L^2$
being $(2p-n)^2$ it follows that the only degrees present in $\z$ are $0$ and $n$.  Therefore $\z$ is a linear combination of $1$ and $\nu$.
\end{proof}

\begin{lema} \label{trace}Let $\SS$ be any irreducible, real, $\Cl_n$ module where $n \equiv 0,3$ (mod 4). For any $\phi$ in $\Cl_n$ such that $\alpha(\phi^t)=\phi$ we have \\
(i) $Tr(\mu_{\phi}) = dim_{\mathbb{R}} \SS \, \langle \phi, 1 \rangle$, if $n \equiv 0$ (mod 4), \\
(ii) $Tr(\mu_{\phi}) = 2 \, dim_{\mathbb{R}} \SS \, \langle \phi, 1 \rangle$, if $n \equiv 3$ (mod 4) .
\end{lema}
\begin{proof}
Let us consider the linear sub-space of $\Cl_n$ given by $\mathcal{S} = \{\phi : \alpha(\phi^t) = \phi \}$. We now pick some orthonormal basis $\{e_i : 1 \le i \le n \}$ in $V$ and
observe that $e^{ij} \mathcal{S} e^{ij} \subseteq \mathcal{S}$, where $e^{ij} = e^i \wedge e^j, i \neq j$. Let now $t :\mathcal{S} \to \mathbb{R}$ be given as $t(\phi) =
Tr(\mu_{\phi})$ for all $\phi$ in $\mathcal{S}$. Since this is linear, it can be written as $t = \langle \cdot, T \rangle$ for some $T$ in $\mathcal{S}$. From the independence of the
trace of some orthonormal basis in $\SS$ and \eqref{spinorscalarprod} we get $t(e^{ij}\phi e^{ij}) = -t(\phi)$ for all $\phi$ in $\mathcal{S}$. Using \eqref{scprod} this results in
having $e^{ij}Te^{ij} = -T$ or further $[T,e^{ij}] = 0$ for all $1 \le i \neq j \le n$, where we have used that $(e^{ij})^2 = -1$. Henceforth $[T, \Lambda^2] = 0$, leading by Lemma
\ref{prep} to $T = \lambda_1 + \lambda_2 \nu$ for some $\lambda_1, \lambda_2$ in $\mathbb{R}$ which can be computed as $\lambda_1 = Tr(\mu_{1}) = dim_{\mathbb{R}}\SS$ and $\lambda_2 =
Tr(\mu_{\nu}) = dim_{\mathbb{R}} {\SS}^{+} - dim_{\mathbb{R}} \SS^{-} = 0$ for $n \equiv 0$ (mod 4) and accordingly $Tr(\mu_{\nu}) = dim_{\mathbb{R}} \SS$ for $n \equiv 3$ (mod 4),
where we have used that $\nu$ belongs to $\mathcal{S}$.
\end{proof}

\section{Structure results}

\subsection{The general setup}

In this section our aim is mainly to locate some classes of holonomy algebras inside the Clifford algebra and derive a number of general properties they must satisfy. Let $A$ be the
subset of $\Cl_n$ given by
$$ A=\{ \phi \in \Cl_n : \phi^t = -\phi\}.$$
This is meant to be the model algebra for most classes of holonomy algebras we will be looking at, in a sense to be made precise below.

\begin{lema} \label{strA} The following hold :
\begin{itemize}
    \item[(i)]  $A$ is a Lie sub-algebra of $(\Cl_n, [\cdot, \cdot])$.
    \item[(ii)] The symmetric bilinear form $\beta(\phi_1, \phi_2)=\langle \phi_1, \alpha(\phi_2) \rangle$ is non-degenerate on $A$ and invariant, that is
                $$ \beta([\phi_1, \phi_2], \phi_3) = -\beta([\phi_1, \phi_3], \phi_2) $$
                whenever $\phi_k ,1 \le k \le 3$ belong to $A$.
\end{itemize}
\end{lema}

\begin{proof}
(i) Follows immediately from anti-symmetrising that $(\phi_1 \phi_2)^t=\phi_2^t \phi_1^t=\phi_2 \phi_1$ whenever $\phi_1, \phi_2$ belong to
$A$. \\
(ii) The non-degeneracy of $\beta$ follows from $A$ being preserved by the involution $\alpha$. Now
\begin{equation*}
\begin{split}
\beta([\phi_1, \phi_2], \phi_3)=&\langle [\phi_1, \phi_2], \alpha(\phi_3) \rangle=\langle \phi_1 \phi_2-\phi_2 \phi_1 , \alpha(\phi_3)\rangle\\
=&\langle \phi_2, \alpha(\phi_1^t) \alpha(\phi_3)-\alpha(\phi_3)\alpha(\phi_1^t)\rangle\\
=&\langle \phi_2, \alpha([\phi_1, \phi_3])\rangle=-\beta(\phi_2, [\phi_1,\phi_3]).
\end{split}
\end{equation*}
\end{proof}

Since $A$ is stable under $\alpha$, it inherits from $\Cl_n$ a bi-grading $A = A^0 \oplus A^1$ into its even respectively its odd degree components. The usual rules $[A^0, A^0]
\subseteq A^0, \ [A^0, A^1] \subseteq A^1, \ [A^1,A^1] \subseteq A^0$ apply, in particular $A^0$ is a Lie sub-algebra of $A$. We can obtain now first order information about some of
the holonomy algebras, by assuming the generating element to be well related to the standard decompositions of $\Cl_n$.

\begin{pro} \label{loc1}
Let $T$ belong to $\Cl_n^0$ and satisfy $T^t = T$. Then :
\begin{itemize}
    \item[(i)]  $\mathfrak{g}_{T}^{\star}$ is a Lie sub-algebra of $A$
    \item[(ii)] $\alpha(\mathfrak{g}_{T}^{\star})=\mathfrak{g}_{T}^{\star}$.
\end{itemize}
\end{pro}
\begin{proof} (i) follows eventually after checking that the generating set $\{X \lrcorner T : X \in V \}$ is contained in $A$ as
$$ (X \lrcorner T)^t=-X \lrcorner T ^t=-X \lrcorner T$$
for all $X$ in $V$. To prove (ii) we notice that $\alpha(X \lrcorner T)=-X \lrcorner T$ for all
$X$ in $V$, in other words $\alpha$ preserves the generating set. Since $\alpha$ is a Lie algebra automorphism preserving $A$, the result follows.
\end{proof}
Therefore, for any $T$ in $Cl_n$, the Lie algebra $\mathfrak{g}_T^{\star}$ splits as
$$\mathfrak{g}_T^{\star}=\mathfrak{g}_T^{\star,0} \oplus \mathfrak{g}_T^{\star,1},$$
where the obvious notation applies. Getting closer to the specific features of the algebras $\mathfrak{g}_T^{\star}$ requires some Lie algebra background we shall now briefly outline. For our setup most convenient is
to adopt the following
\begin{defi} Let $\mathfrak{g}$ be a real Lie algebra. It is called semisimple if it admits a symmetric bilinear form $\beta$ which is non-degenerate and satisfies
\begin{equation*}
\beta([\phi_1, \phi_2], \phi_3) = - \beta([\phi_1, \phi_3], \phi_2)
\end{equation*}
for all $\phi_1, \phi_2, \phi_3$ in $\mathfrak{g}$.
\end{defi}
This essentially ensures that any ideal $\mathfrak{i}$ of $\mathfrak{g}$ has trivial extension, that is there exists $\mathfrak{i}^{\perp}$ such that $\mathfrak{g}= \mathfrak{i} \oplus
\mathfrak{i}^{\perp}$, where $\mathfrak{i}^{\perp}$ denotes the orthogonal complement of $\mathfrak{i}$ w.r.t the non-degenerate form $\beta$. In particular
\begin{pro} \label{crit}
Let $\mathfrak{g} $ be a real Lie algebra. If $\mathfrak{g}$ is semisimple then $[\mathfrak{g}, \mathfrak{g}] = \mathfrak{g}$ if and only if it has trivial center. Here the center
$Z(\mathfrak{g})$ of $\mathfrak{g}$ is given as $Z(\mathfrak{g}) = \{ \z \in \mathfrak{g} : [\z,\mathfrak{g}] = 0 \}$.
\end{pro}
Real Lie algebras $\mathfrak{g}$ satisfying $[\mathfrak{g}, \mathfrak{g}]=\mathfrak{g}$ are termed {\it{perfect}} and Proposition \ref{crit} provides a criterion for checking this, to be used later on.
\begin{lema} \label{propA}
Suppose that $n \equiv 0$ (mod 8). The following hold:
\begin{itemize}
    \item[(i)]      $A$ is a semi-simple Lie algebra.
    \item[(ii)]     $A$ isomorphic to $\mathfrak{so}(d,d)$ where $d=\frac{1}{2}dim_{\mathbb{R}}\SS $ and $\SS$ is the irreducible real $\Cl_n$ module.
    \item[(iii)]    the adjoint representation of $A^0$ on $A^1$ is irreducible.
\end{itemize}
\end{lema}

\begin{proof} (i) follows from the non-degeneracy of $\beta$ on $A$, which is due to $\alpha(A)=A$. \\
(ii) Let us equip $\SS$ with the scalar product $\hat{\beta}$ which leaves $\SS^{+}$ and $\SS^{-}$ orthogonal and equals $\pm \langle \cdot, \cdot \rangle$ on $\SS^{\pm}$. In short,
$\hat{\beta}(x,y)= \langle \nu x, y \rangle $ for all $x, y$ in $\SS$. If $\phi$ is in $A$ it is easy to check that $\hat{\beta}(\mu_{\phi}x, y) + \hat{\beta}(\mu_{\phi}y, x) = 0$ for
all $x, y$ in $\SS$, that is $\mu_{\phi}$ belongs to $\mathfrak{so}(\SS, \hat{\beta}) \cong \mathfrak{so}(d,d)$. But when $n \equiv 0$ (mod 8) the Clifford multiplication gives a
linear isomorphism $\mu: \Cl_n \to End(\SS,\SS)$ which is also a Lie algebra isomorphism and our claim follows. \\
(iii) follows standardly from (ii).
\end{proof}
Similar results can be proved in the remaining series of dimensions but this is somewhat beyond the scope of the present paper. In the same vein
\begin{pro} For any $T$ in $\Cl_n^0$ with $T^t = T$ the Lie algebra $\mathfrak{g}_T^{\star}$ is semisimple.
\end{pro}
\begin{proof} We need only see that the restriction of $\beta$ to $\mathfrak{g}_T^{\star}$ is non-degenerate. But this follows easily from the fact that $\alpha$ preserves
$\mathfrak{g}_T^{\star}$.
\end{proof}
\begin{lema} \label{flip}
Let $T$ be in $\Cl_n^0$. Then
\begin{equation*}
    \mathfrak{g}_{eTe}^{\star} = e \, \mathfrak{g}_{T}^{\star} \, e
\end{equation*}
for any unit vector $e$ in $V$.
\end{lema}

\begin{proof}
At first we notice that $eTe$ still belongs to $\Cl_n^0$. We have
\begin{equation*}
-2 X \lrcorner (eTe) = X(eTe)-(eTe)X = -e(XT)e+e(TX)e = 2e(X \lrcorner T)e
\end{equation*}
for all $X$ in $(e)^{\perp}$. Similarly, $e \lrcorner (eTe)=e(e \lrcorner T)e$ hence $X \lrcorner (eTe)=e(F_eX \lrcorner T)e$ for all $X$ in $V$, where $F_e$ is the invertible
endomorphism of $V$ which equals $-1$ on $(e)^{\perp}$ and $1$ on $(e)$. Let $\rho_e: \Cl_n \to \Cl_n$ be defined as $\rho_e(\phi)=e \phi e$ for all $\phi$ in $\Cl_n^0$. It therefore
maps the generating set of $\mathfrak{g}_T^{\star}$ onto that of $e\mathfrak{g}_T^{\star}e$ and since $-\rho_e$ is a Lie algebra isomorphism it is easy to conclude.
\end{proof}

When $n \equiv 0$ (mod 4) the map $\rho_e$ intertwines, say, $\Cl_n^{+} \cap \Cl_n^0$ and $\Cl_n^{-} \cap \Cl_n^0$, therefore from the Lemma above we see that the holonomy algebra does
not distinguish between the generating elements being in $\Cl_n^{+}$ or $\Cl_n^{-}$. Hence all results obtained for holonomy algebras generated by elements in $\Cl_n^{+}$ extend
automatically to generating elements in $\Cl_n^{-}$. We end this section by an example of forms when the holonomy algebras can be easily computed.

\begin{pro} \label{volforms}
Let $(V^n, \langle \cdot , \cdot \rangle)$ be an Euclidean vector space oriented by $\nu$ in $\Lambda^n(V)$. Then:
\begin{itemize}
    \item[(i)]      $\g^{\star}_\nu = \solie(n, 1)$ for $n \equiv 0, 1$ (mod 4)
    \item[(ii)]     $\g^*_\nu = \solie(n+1)$ for $n \equiv  2, 3$ (mod 4)
    \item[(iii)]    in both cases $\g^*_\nu \psi = 0$ if and only if $\psi = 0$, for any $\psi \in \SS$.
\end{itemize}
\end{pro}
\begin{proof}
Let us first notice that the generating set $\{X \hook \V:X \in V\}$ is isomorphic to $V$, since the volume form $\V$ is non-degenerate.
Keeping in mind that by \eqref{clifford1form} we have $X \hook \V = - X \nu$ and
using \eqref{vsquare}, \eqref{alphaeven} this yields
    $$ [X \hook \nu, Y \hook \nu] = (-1)^{\frac{1}{2}(n+1)(n+2)} [X, Y]=2 (-1)^{\frac{1}{2}(n+1)(n+2)} X\wedge Y$$
for all $X,Y$ in $V$. Similarly we get for the triple commutators
    $$ [\alpha, X \hook \nu] = [X, \alpha] \nu = 2 \, FX \hook \nu $$
for all $X$ in $V$, where $\alpha=\langle F \cdot, \cdot \rangle$ belongs to $\Lambda^2(V)$. Therefore $\g^*_\nu= V \oplus \Lambda^2(V)$ as a vector space and the claims in (i) and
(ii) follow from the commutator rules above. (iii) follows easily from the invertibility of $\V$ in $Cl_n$, as defined in Def. \ref{invert}.
\end{proof}

\subsection{The set of fixed spinors}
As it will appear below the holonomy algebra of some element $T$ in $\Cl_n$ is intimately related to the space of spinors fixed by $T$, which we recall
to be defined as
\begin{equation} \label{spininv}
Z_T=\{\psi \in \SS: (X \lrcorner T) \psi = 0 , \  \mbox{for all} \ X \in V \}.
\end{equation}
Notice that if (the non-zero) $T$ is of degree $1$ or $2$ the set $Z_T$ is obviously reduced to zero and moreover the latter holds for forms of degree $3$ (see \cite{holo}). We now
gather a number of basic facts concerning the set $Z_T$. If $n \equiv 0$ (mod 4) we split $Z_T$ along the splitting $\SS = \SS^{+} \oplus \SS^{-}$ and get
\begin{equation*}
Z_T=Z_{T}^{+} \oplus Z_T^{-}
\end{equation*}
where the obvious notation applies.

\begin{lema} \label{z1}
Let $T$ belong to $\Cl_n^0 \cap \Cl_n^{+}$ where $n \equiv 0$ (mod 4). Then
\begin{itemize}
\item[(i)]  \begin{equation*}
                \begin{split}
                Z_T^{+} &= \{\psi \in \SS^{+} : T \psi = 0\} \\
                Z_T^{-} &= \{\psi \in \SS^{-} : T V \psi = 0\}
                \end{split}
            \end{equation*}
\item[(ii)] if $n=8$, then $Z_T^{-}=(0)$ provided $T$ does not vanish.
\end{itemize}
\end{lema}

\begin{proof}
(i) follows directly from the stability conditions. \\
(ii) If $Z_{T}^{-} \neq (0)$ there exists a non-zero $\psi$ in $\SS^{-}$ with $T V\psi=0$. But in $8$-dimensions $V\psi=\SS^{+}$ hence $T \SS^{+}=0$ whence $T=0$, a contradiction.
\end{proof}
\begin{defi} \label{invert}An element $T$ of $\Cl_n^0 \cap \Cl_n^{+}$ is invertible in $\Cl_n \slash (1,\nu)$ if there exists $T^{-1}$ in $\Cl_n^0 \cap \Cl_n^{+}$ with $T^{-1}T = T T^{-1} = 1 + \nu$.
\end{defi}
Actually a necessary and sufficient condition for $Z_T$ to vanish is
\begin{pro} \label{z2}
Let $T$ belong to $\Cl_n^0 \cap \Cl_n^{+}, n \equiv 0$ (mod 4) such that $T^t = T$. Then $Z_T = (0)$ if and only if $T$ is invertible in $\Cl_n  \slash (1,\nu)$.
\end{pro}
\begin{proof} If $T$ is invertible Lemma \ref{z1} yields immediately the vanishing of $Z_T^{\pm}$ hence that of $T$.  Suppose now that $Z_T=(0)$ and let $L_T: \Cl_n^{+} \cap \Cl_n^0
\to \Cl_n^{+} \cap \Cl_n^0$ be left multiplication with $T$. If $\phi $ is in the kernel of $L_T$ it follows that $T( \phi {\SS^{+}} ) = 0$ and moreover, since $\phi \SS^{+} \subseteq
\SS^{+}$ Lemma \ref{z1} tells us that $\phi \SS^{+} \subseteq Z_T$ hence $\phi {\SS^{+}}=0$. Therefore $\phi $ vanishes and it follows that $L_T$ is injective, thus invertible and this
provides easily the required inverse for $T$ in $\Cl_n \slash (1, \nu) $ given that $1_{\Cl_n^{+} \cap \Cl_n^0} = \frac{1}{2} L_{1 + \nu}$.
\end{proof}

In the rest of this section we shall present examples of situations when the set of fixed spinor can be seen directly to be trivial.

\begin{pro} \label{2sq}
Let $\alpha$ in $\Lambda^2(V)$ be a two-form such that $T = \alpha \wedge \alpha \neq 0$. Then $Z_T = (0)$.
\end{pro}

\begin{proof}
Let $F$ be the skew-symmetric endomorphism associated to $\alpha$ via the metric $g$, that is $\alpha = \langle F \cdot, \cdot \rangle$. We have $X \hook (\alpha \wedge \alpha) = 2 (X
\hook \alpha) \wedge \alpha = 2FX \wedge \alpha$ for all $X$ in $V$. Let now $\psi$ be in $Z_T$ and let us set  $r = {\rm rank}(F)$. From
\begin{equation*}
    (X \hook T) \psi = 2 \; (FX \wedge \alpha) \psi = 0
\end{equation*}
follows $(X \wedge \alpha) \psi = 0$ for all $X$ in  $Im(F)$. By Clifford contraction we get immediately
    $$ \sum_{e_i \in \, Im (F)} e_i (e_i \wedge \alpha) \psi = -\sum_{e_i \in \, Im (F)} (e_i \hook (e_i \wedge \alpha)) \psi = (2 - r) \; \alpha \psi = 0, $$
which leads to $\alpha \psi = 0$ since having $r = 2$ would imply $T = \alpha \wedge \alpha = 0$, a contradiction. Therefore
    $$ 0 = X \alpha \psi = (X \wedge \alpha - X \hook \alpha) \psi $$
for all $X$ in $V$ and since $(X \wedge \alpha) \psi = 0$ for all $X$ in $Im(F)$ we are lead to $(X \lrcorner \alpha)\psi = 0$ for all $X$ in $Im(F)$. It follows that $\psi = 0$ as
$\alpha \neq 0$.
\end{proof}

A very simple observation, which appears to be useful in low dimensions is

\begin{lema} \label{trickl}
Let $T$ belong to $\Lambda^k(V), k \neq 0$. Then: \\
(i)  $T Z_T = 0$\\
(ii) $Z_T = Z_{\nu T}$
\end{lema}

\begin{proof}(i) If $\psi$ belongs to $Z_T$ we have $(X \lrcorner T)\psi = 0$ for all $X$ in $V$. Therefore
$ \sum \limits_{i = 1}^{n} e_i (e_i \lrcorner T) \psi = 0$
for some orthonormal frame $\{e_i, 1 \le i \le n\}$ leading to $k T \psi = 0$ and the claim follows.\\
(ii) From (i) we get that $Z_T = \{ \psi \in \SS : T V \psi = 0 \}$ for any pure degree form $T$ and the claim follows easily.
\end{proof}

\begin{pro} \label{ldim}
For any $T$ in $\Lambda^k(\mathbb{R}^n), n \leq 7$  the set $Z_T$ is trivial.
\end{pro}

\begin{proof} This is obvious when $k=1,2$ and when $k=3$ it was proved in \cite{holo}. Now if $k \ge 4$ we have $Z_T = Z_{\nu T}$ and since $\nu T$ has degree
$n-k \le 3$ we conclude by the above.
\end{proof}

Therefore the first case of interest is that of dimension $8$, which will be studied in detail in the latter part of the paper.

\subsection{Perfect fix algebras}

In this section we shall examine situations when the fix algebra $\mathfrak{g}_T^{\star}$ for some $T$ in $\Cl_n$ is perfect. We will see that this is not always the case since those
can be abelian by the examples in the next section. However we will show that it is possible to give necessary conditions to that extent. Let us set first a preparatory Lemma.

\begin{lema} \label{prepc}
If $\z_1, \z_2$ in $\Cl_n$ satisfy $\z_1 X=X \z_2$ for all $X$ in $V$ then $\z_1$ and $\z_2$ belong to $(1, \nu)$.
\end{lema}

\begin{proof}
It follows that $-\vert X \vert^2 \zeta_1=X \z_2 X$ for all $X$ in $V$ and further
\begin{equation*}
\vert X \vert^2 Y\z_2 Y=\vert Y\vert^2 X\z_2 X
\end{equation*}
for all $X,Y$ in $V$. By left multiplication with some non-zero $X$ we get $(XY)\z_2 Y=-\vert Y \vert^2 \z_2 X$ and now right multiplication with a non-zero $Y$ yields
$(XY)\z_2=\z_2 (XY)$ for all $X,Y$ in $V$.  Now $\z_2$ is in $(1, \nu)$ by Lemma \ref{prep} and it is easy
to see this implies the claim for $\z_1$ as well.
\end{proof}

\begin{pro} \label{center}
Suppose that $n \equiv 0$ (mod 4) and let $T$ satisfy $\nu T = T \nu = T$ and $T^t = T$.  If $Z_T = (0)$ then $\mathfrak{g}_{T}^{\star}$ has trivial center.
\end{pro}

\begin{proof} If $\z$ in $Z(\mathfrak{g}_T^{\star})$ we must clearly have
\begin{equation} \label{centm}
[\z, XT-TX]=0
\end{equation}
for all $X$ in $V$. Since $X T-TX$ belongs to $\Cl_n^1$ for all $X$ in $V$ by applying $\alpha$ to the equation above we get that $[\alpha(\z), XT-T X]=0$ for all $X$ in $V$, hence
after splitting $\z$ into its even resp. odd components it is enough to treat \eqref{centm} when $\z$ belongs to $\Cl_n^0$ resp. $\Cl_n^1$.\\
{\bf{Case I}}: $\z$ belongs to $\Cl_n^0$. \\
Since $\nu(XT-TX)=-(XT+TX)$ for all $X$ in $V$ and $\nu \z=\z \nu$ after left multiplication of \eqref{centm} with the volume form we get $\z(XT+TX)=(XT+TX)\z$ whenever $X$ belongs to $V$.
Taking linear combinations with \eqref{centm} gives further
\begin{equation*}
\begin{split}
\z XT =&X T\z \\
\z TX=& TX \z
\end{split}
\end{equation*}
for all $X$ in $V$. Since $Z_{T}=(0)$ we know that $T$ must be invertible in $\Cl_n \slash (1, \nu)$ (see Proposition \ref{z2}) hence using the second equation above we have
\begin{equation*}
(T^{-1}\z T)X=T^{-1}TX\z=(1+\nu)X\z=X(1-\nu)\z
\end{equation*}
for all $X$ in $V$. But from Lemma \ref{prepc}, (ii) we get that $(1-\nu)\z$ belongs to $(1,\nu)$ and hence vanishes as $((1-\nu)\z)^t=-(1-\nu)\z$. From the vanishing of $(1-\nu)\z$ it follows that
$T^{-1}\z T=0$ and this leads after right multiplication with $T^{-1}$ resp. left multiplication with $T$ to $(1+\nu)\z=0$. Thus $\z=0$ in this case. \\
{\bf{Case II}}: $\z$ belongs to $\Cl_n^1$. \\
We have as before $\z(XT-TX)=(XT-TX)\z$ for all $X$ in $V$.  But $\nu (XT-TX)=-(XT+TX)$ for all $X$ in $V$ because $\nu T=T$ and since $\nu \z =-\z \nu$ after left multiplication with the volume
form we obtain
$$ \z (XT+TX)=-(XT+TX) \z
$$
for all $X$ in $V$. Taking linear combinations with the original equation gives now
\begin{equation*}
\begin{split}
(\z T)X=&-X (T \z)\\
\z XT=&-TX \z
\end{split}
\end{equation*}
for all $X$ in $V$.  Using Lemma \ref{prepc}, (ii) we then get that $\z T$ and $T \z$ belong to $(1,\nu)$ and therefore must vanish since elements of $\Cl_n^1$.  The invertibility of
$T$ in $\Cl_n \slash (1,\nu)$ leads then to $(1+\nu)\z=0$ whence left multiplication by $1+\nu$ in the second equation above gives $TX\z=0$ for all $X$ in $V$. Again the invertibility
of $T$ implies that $(1+\nu)X\z=X(1-\nu)\z=0$ for all $X$ in $V$ and we conclude that $(1-\nu)\z=0$  hence $\z=0$ and the proof is finished.
\end{proof}

Summarising, after making use of the semisimplicity of $\mathfrak{g}_T^{\star}$ and of Proposition \ref{crit}, we obtain

\begin{teo}\label{perf1} Let $T$ belong to $\Cl_n^0 \cap \Cl_n^{+}$ where $n \equiv 0$ (mod 4) and satisfy $T^t = T$. If moreover $Z_T = (0)$, the algebra $\mathfrak{g}_T^{\star}$ is perfect, that is
\begin{equation*}
\mathfrak{g}_T^{\star} = [\mathfrak{g}_T^{\star}, \mathfrak{g}_T^{\star}].
\end{equation*}
\end{teo}

\section{Holonomy algebras from distinguished Clifford algebra elements}

\subsection{Unipotent elements}

In this section we shall compute directly the fix and holonomy algebras of a unipotent element $T$ of $Cl_n^{+}, n \equiv 0$ (mod 4) as introduced below.

\begin{defi} \label{unipotdef}
Let $T$ belong to $\Cl_n^{+}$ where $n \equiv 0$ (mod 4). It is called unipotent if it satisfies $T^t = T$ and $T^2 = 1 + \nu$.
\end{defi}
In particular any unipotent element $T$ belongs to $Cl_n^0$. We need first to state and prove the following preliminary result, to be used later on as well.
\begin{lema} \label{dcom}
Let $T$ belong to $\Cl_n^{+} \cap \Cl_n^0$ where $n \equiv 0$ (mod 4). Then :
\begin{equation*}
4[X \lrcorner T, Y \lrcorner T]=-T[X,Y]T+YT^2X-XT^2Y
\end{equation*}
whenever $X,Y$ belong to $V$.
\end{lema}

\begin{proof}
Follows directly from the stability relations under the form $TXT=0$ for all $X$ in $V$. Details are left to the reader.
\end{proof}

\begin{teo} \label{uni} Let $T$ be a unipotent element $T$ of $\Cl_n^{+}$ where $n \equiv 0$ (mod 4). Then:
\begin{itemize}
\item[(i)] $Z_T=(0)$
\item[(ii)] the fix algebra of $T$ is perfect
\item[(iii)] the holonomy algebra of $T$ is isomorphic to $\mathfrak{so}(n,1)$.
\end{itemize}
\end{teo}

\begin{proof}
(i) Since any unipotent element $T$ of $\Cl_n^{+}$ is clearly invertible in $\Cl_n \slash (1,\nu)$ Proposition
\ref{z2} implies that $Z_T = (0)$.\\
(ii) follows from (i) and Theorem \ref{perf1}. \\
(iii) For notational convenience let $E_T = \{X \lrcorner T: X \in V\}$ be the generating set of $\mathfrak{g}_T^{\star}$. It is isomorphic to $V$ under the map $\iota^1 : V \to E_T,
\iota^1(X) = X \lrcorner T$. Here only the injectivity of $\iota^1$ has to be proved, and indeed, if $\iota^1(X) = 0$ it follows that $XT = TX$ and further $0 = TXT = T^2X$ leading to
the vanishing of $X$. Now by Lemma \ref{dcom} combined with the unipotency of $T$ the space $[E_T, E_T]$ equals $\{ T\alpha T + (1 - \nu)\alpha: \alpha \in \Lambda^2(V)\}$. This is
isomorphic to $\Lambda^2(V)$ under $\iota^2 : \Lambda^2(V) \to [E_T, E_T], \iota^2(\alpha) = T \alpha T + (1 - \nu) \alpha$. Indeed, if $\iota^2(\alpha) = 0$ we find that $T \alpha T +
(1 - \nu) \alpha = 0$ but then both summands vanish as the first is in $\Cl_n^{+}$ and the second in $\Cl_n^{-}$. Hence $\alpha = 0$ and so $\iota^2$ is injective. Moreover
\begin{equation*}
\begin{split}
[T\alpha T+(1-\nu)\alpha, T\beta T+(1-\nu)\beta]=&[T\alpha T, T\beta T]+[(1-\nu)\alpha, (1-\nu)\beta]\\
+&[(1-\nu)\alpha, T\beta T]+[T\alpha T, (1-\nu)\beta]\\
=&[T\alpha T, T\beta T]+[(1-\nu)\alpha, (1-\nu)\beta ]
\end{split}
\end{equation*}
for all $\alpha, \beta$ in $\Lambda^2(V)$ after using that $T$ belongs to $\Cl_n^{+}$. Obviously
$$[(1-\nu)\alpha, (1-\nu)\beta ]=(1-\nu)^2[\alpha, \beta]=2(1-\nu)[\alpha, \beta]$$ and
moreover the unipotency of $T$ leads easily to $[T\alpha T, T\beta T]=(1+\nu)T[\alpha, \beta]T=2T[\alpha, \beta]T$. Altogether this yields
\begin{equation*}
2\iota^2[\alpha, \beta]=[\iota^2\alpha, \iota^2 \beta]
\end{equation*}
for all $\alpha, \beta$ in $\Lambda^2(V)$, in other words $\frac{1}{2}\iota^2 : \Lambda^2(V) \to [E_T, E_T]$ is a Lie algebra isomorphism. Now, we compute
\begin{equation*}
\begin{split}
-2[T\alpha T+(1-\nu)\alpha, X \lrcorner T]=&[T\alpha T+(1-\nu)\alpha, XT-TX]\\
=& [T\alpha T, XT]+[(1-\nu)\alpha, XT]\\&-[T\alpha T, TX]-[(1-\nu)\alpha, TX]
\end{split}
\end{equation*}
We now estimate  each term separately. We have
$$ [T \alpha T, XT]=-XT^2\alpha T=-X(1+\nu) \alpha T= -2X\alpha T
$$
after using $TXT=0$. Similarly, $[T\alpha T, TX]=-2T\alpha X$ and using furthermore that $T$ belongs to $\Cl_n^{+}$ we finally obtain
\begin{equation*}
-2[T\alpha T+(1-\nu)\alpha, X \lrcorner T]=-2X\alpha T+2\alpha X T-2T\alpha X+2TX \alpha.
\end{equation*}
This ends by saying that $[\iota^2\alpha, \iota^1X]=\iota^1[X, \alpha]$ whenever $\alpha$ belongs to $\Lambda^2(V)$ and $X$ in $V$. Summarising, it follows that
$\mathfrak{g}_T^{\star}=E_T \oplus [E_T, E_T]$ and moreover
$$\iota^1 \oplus \frac{1}{2} \iota^2 : V \oplus \Lambda^2(V) \to \mathfrak{g}_T^{\star}$$
realises the desired Lie algebra isomorphism with $\mathfrak{so}(n,1)$.
\end{proof}

\begin{rema}
(i) We shall see in the next section that unipotent elements naturally play a special role in the classification of holonomy algebras of self-dual $4$-forms in dimension $8$.\\
(ii) Explicit examples of unipotent elements are easy to make. When the dimension of our vector space $V$ satisfies $dim_{\mathbb{R}}V \equiv 0$ (mod 4) we see that
$\frac{1}{\sqrt{2}}(1 + \nu)$ is unipotent and therefore Theorem \ref{uni} recovers partly results in Proposition \ref{volforms}. Moreover, if we take $V_1,V_2$ to be Euclidean vector
spaces of dimensions $\equiv 0$ (mod 4) oriented by volume forms $\nu_k, k=1,2$. Then $\frac{1}{\sqrt{2}}(\nu_1+\nu_2)$ is an unipotent element of the direct product space $V_1 \times V_2$.

\end{rema}

\subsection{Squares of spinors}

We shall first recall in what follows some facts about the squaring construction in two series of dimensions. To begin with, let $(V^n, \langle \cdot, \cdot \rangle)$ be a Euclidean
vector space which furthermore is supposed to be oriented, with orientation given by $\nu$ in $\Lambda^n(V)$. A peculiar property of the Clifford multiplication when $n \equiv 8$ (mod 8) is then to give
an isomorphism (see \cite{laws}):
\begin{equation} \label{cliffordiso}
    \mu: \Cl_n \to {\rm Hom}_\Real(\SS, \SS)
\end{equation}
where $\SS$ is the irreducible real $Cl_n$ module. When $n \equiv 7$ (mod 8) this still holds provided $Cl_n$ is replaced by $Cl_n^{+}$. Let us now fix a spinor $x \in \SS^+$ (or in
$\SS$ if $n \equiv 7$ (mod 8)), which we normalise to $|x|= 1$. Then the isomorphism \eqref{cliffordiso} gives rise to an element $x \otimes x \in \Cl_n$ (or $Cl_n^{+}$ when $n \equiv
7$ (mod 8)) such that:
\begin{equation} \label{spinorsquare}
    (x \otimes x) \, \psi = \langle \psi, x \rangle \, x
\end{equation}
for all $\psi$ in $\SS$. The element $x \otimes x$ consists of forms of various degrees and is customarily called the square of $x$. Indeed it is well known \cite{laws} that
\begin{equation} \label{spinorsquarecontent}
    x \otimes x = \sum \limits_{k \equiv 0,3 \ ({\rm mod} \ 4)}^{n} (x \otimes x)_k,
\end{equation}
where $(x \otimes x)_k$ denotes the projection of $(x \otimes x)$ onto $\Lambda^{k}(V)$. Note that when $n \equiv 0$ (mod 8) the odd degrees are not present. Below we list
some of the properties of $x \otimes x$, of relevance for our study.
\begin{lema} \label{squareproperties}
    Let $x$ be a unit length spinor in $\SS^+$, where
    $n \equiv 7,8$ (mod $8$). The following hold:
    \begin{itemize}
    \item[(i)]      the spinor square $x \otimes x$ is an idempotent of $\Cl_n$, that is $(x \otimes x)^2=x \otimes x$.
    \item[(ii)]     if $n \equiv 0$ (mod 8) we have $\nu \, (x \otimes x) = (x \otimes x) \, \nu = x \otimes x$.
    \item[(iii)]    for all $\vph \in \Cl_n$ we have
                    $$ (x \otimes x) \, \vph \, (x \otimes x) = \kappa \, \langle \vph, x \otimes x \rangle \, (x \otimes x), $$
                    where $\kappa = 2^{\frac{n}{2}}$ in $n = 0$ (mod 8) and $\kappa = 2^{\frac{n+1}{2}}$ when $n \equiv 7$ (mod 8).
    \end{itemize}
\end{lema}
\begin{proof}
(i) We use \eqref{spinorsquare} for $\psi = x$ which gives $(x \otimes x) \, x = x$. Therefore left multiplication of \eqref{spinorsquare} with $x \otimes x$ gives:
$$ (x \otimes x)^2 \, \psi = \langle \psi, x \rangle \, (x \otimes x) \, x = \langle \psi, x \rangle \, x = (x \otimes x) \, \psi, $$
for all $\psi$ in $\SS$ and the claim follows.\\
(ii) We use the definition \eqref{spinorsquare} to obtain after recalling that $x \in \SS^+$
$$ \nu \, (x \otimes x) \, \psi = \langle \psi, x \rangle \, \nu \, x = \langle \psi, x \rangle \, x = (x \otimes x) \, \psi, $$
for all $\psi \in \SS$. Since $x \otimes x \in \Cl^0_n$ and $n \equiv 0$ (mod 8), we further have $[\nu, x \otimes x] = 0$.\\
(iii) Let $\phi$ belong to $\Cl_n$ and $n \equiv 0$ (mod 8). Using again \eqref{spinorsquare} we compute
\begin{equation*}
\begin{split}
(x \otimes x) \, \vph \, (x \otimes x) \, \psi=& \langle \psi, x \rangle \, (x \otimes x) \, \vph \, x =\langle \psi, x\rangle \langle \phi x,x \rangle x\\
=& \langle \phi x,x \rangle (x \otimes x)\psi
\end{split}
\end{equation*}
for all $\psi$ in $\SS$ hence $(x \otimes x)\phi (x \otimes x)= \langle \phi x,x \rangle x \otimes x$ for all $\phi$ in $\Cl_n$. Recall now \cite{dado} that $\langle 1, x \otimes x
\rangle = 2^{-\frac{n}{2}}$ and , fact which follows essentially by taking traces and using Lemma \ref{trace}. Therefore
\begin{equation*}
\begin{split}
2^{-\frac{n}{2}}\langle \phi x,x \rangle=& \langle 1, (x \otimes x)\phi (x \otimes x)\rangle=\langle \alpha(x \otimes x)^t \cdot 1, \phi (x \otimes x)\rangle \\
=& \langle x \otimes x, \phi (x \otimes x)\rangle=\langle (x \otimes x) \alpha(x \otimes x)^t, \phi \rangle\\
=& \langle (x \otimes x)^2, \phi \rangle= \langle x \otimes x, \phi \rangle
\end{split}
\end{equation*}
and the claim follows. For $n \equiv 7$ (mod 8) this is proved analogously.
\end{proof}

Based on the technical Lemma above we shall compute now the holonomy algebras of the square of a spinor. Let us begin with the case of $n \equiv 0$ (mod 8).

\begin{teo}
Let $n \equiv 0$ (mod 8) and $x$ be a unit length spinor  in $\SS^+$ and $x \otimes x$ be its square. Then:
\begin{itemize}
\item[(i)] $\g^*_{x \otimes x}=\mathfrak{h}_{x \otimes x}^{\star} \cong \solie(n,1)$
\item[(ii)] $Z_{x \otimes x}=\{\psi \in \SS^{+} : \psi \perp x\} \oplus \{ \psi \in \SS^{-} : \psi \perp Vx \}$.
\end{itemize}
\end{teo}
\begin{proof}
(i) Let $E_T=\{X \hook (x \otimes x):X \in V\}$ be the generating set of $V$, which is easily seen to be isomorphic to $V$ under the map $\iota^1 : V \to E_T, \iota^1(X)=X \hook (x
\otimes x)$. Further on,  let us define $\iota^2 : \Lambda^2(V) \to \Cl_n(V)$ by $ \iota^2 (\alpha) = \sum \limits_{i=1}^n  e_i \, (x \otimes x) \, (e_i \hook \alpha)$ for some
orthonormal frame $\{e_i, 1 \le i \le n\}$ on $V$. This is injective since from $\iota^2(\alpha)=0$ we get by right multiplication with $X(x \otimes x), X$ in $V$ to $\sum
\limits_{i=1}^{n} e_i \, (x \otimes x) \, (e_i \hook \alpha) \, X \,  (x \otimes x)=0$ and further by means of Lemma  \ref{squareproperties} to $(X \hook \alpha) x \otimes x=0$ for all
$X$ in $V$, whence $\alpha=0$. Furthermore Lemma \ref{squareproperties} combined with the stability relations leads easily to $[\iota^2\alpha, \iota^2 \beta]=[\alpha, \beta]$ for all
$\alpha, \beta$ in $\Lambda^2(V)$.
\begin{eqnarray*}
        4 \, [X \hook (x \otimes x), Y \hook (x \otimes x)]
        &=& - X (x \otimes x) Y + Y (x \otimes x) X \cr
        & & - (x \otimes x) XY (x \otimes x) + (x \otimes x) XY (x \otimes x) \cr
        &=& - \iota^2(X \wedge Y),
    \end{eqnarray*}
whenever $X, Y$ belong to $V$. Here we have made once more extensive use of the stability conditions and of Lemma \ref{squareproperties},
under the form $(x \otimes x)  \Lambda^2 (x \otimes x)=0$. Hence the even commutators span $\solie(n)$. The
    triple commutator is similarly computed:
    \begin{equation*}
\begin{split}
    [\iota^2(\alpha), X \hook (x \otimes x)]=& -\h \sum \limits_{i=1}^n \Big(e_i (x \otimes x) Fe_i X (x \otimes x) + (x \otimes x) X e_i (x \otimes x) Fe_i \\
    & - X (x \otimes x) e_i (x \otimes x) Fe_i - e_i (x \otimes x) Fe_i (x \otimes x) X \Big) \\
    =& FX \hook (x \otimes x),
\end{split}
    \end{equation*}
where we have made use of Lemma \ref{squareproperties} and have set, for convenience, $\alpha=\langle F \cdot, \cdot \rangle$. Therefore $\iota^1 \oplus \iota^2 : V \oplus \Lambda^2V \to
\mathfrak{g}^{\star}_{x \otimes x}$
gives the desired Lie algebra isomorphism between the fix algebra of $x \otimes x$ and $\mathfrak{so}(n,1)$ which is obviously perfect.\\
(ii) Let $\psi$ belong to $Z_T$. Then by Lemma  \ref{z1} this is equivalent with $(x \otimes x) \psi^{+}=0$ and $(x \otimes x)V \psi^{-}=0$ and the claim follows now from the
definition of $x \otimes x$, where $\psi^\pm \in \SS^\pm$.
\end{proof}

\begin{rema} Squares of spinors provide examples of unipotent elements other than those coming from volume forms.  Indeed, if $x$ belongs to $\SS^{+}$ with $\vert x \vert=1$ it is
easy to see that $2\sqrt{2}(x \otimes x-\frac{1}{4}(1+\nu))$ is a unipotent element in $\Cl_n^{+}$. Despite of the absence of fixed spinors in this case, the holonomy algebra remains
isomorphic to $\mathfrak{so}(n,1)$.
\end{rema}

When $n \equiv 7$ (mod 8) we get fix and holonomy algebras of quite different nature than those seen before. In particular,
those appear not to be perfect.

\begin{teo} \label{abelian}
Let $n \equiv 7$ (mod 8) and $x$ belong to $\SS$ such that $\vert x \vert =1$. Then:
\begin{itemize}
\item[(i)] $\mathfrak{g}_{x \otimes x}^{\star}$ is abelian, isomorphic to $V$ hence $\mathfrak{h}^{\star}_{x \otimes x}=(0)$
\item[(ii)] $Z_{x \otimes x}=(x)^{\perp}$.
\end{itemize}
\end{teo}

\begin{proof}
(i) For any $X,Y$ in $V$ we compute
\begin{equation*}
\begin{split}
4[X \lrcorner (x \otimes x), Y \lrcorner (x \otimes x)]&=(\alpha(x \otimes x)X-
X x \otimes x)(\alpha(x \otimes x)Y-Y x \otimes x)\\
&=\biggl [ \alpha(x \otimes x)X \alpha(x \otimes x) \biggl ]Y-\alpha(x \otimes x)
XY (x \otimes x ) \\
&-X (x \otimes x) \alpha(x \otimes x)Y+X \biggl [ (x \otimes x) Y (x \otimes x) \biggr ].
\end{split}
\end{equation*}
Now $x \otimes x$ belongs to $\Cl_n^{+}$, hence $\alpha(x \otimes x)$ is in $\Cl_n^{-}$ leading to the vanishing of the second and third term above in view of the stability conditions
in Lemma \ref{stability}. Now the first and the last terms vanish too by Lemma \ref{squareproperties}, (iii) and since $x \otimes x$ does not contain degree $1$ forms, therefore
$\mathfrak{g}_{x \otimes x}^{\star}$ is abelian. But $x \otimes x$ is non-degenerate, as it contains a non-zero multiple
of the volume form whence $\mathfrak{g}_{x \otimes x}^{\star}$ is isomorphic with $V$.\\
(ii) follows easily from the construction of $x \otimes x$.
\end{proof}

\section{$8$-dimensions}

In the rest of this paper we shall consider a Euclidean vector space $(V^8, \langle \cdot, \cdot \rangle)$ with orientation given by $\nu$ in $\Lambda^8(V)$. Our aim is to obtain
classification results for holonomy algebras $\mathfrak{h}_T^{\star}$ generated by $T$ in $\Lambda^4_{+}(V)$. Here we recall that in dimension $8$ the Hodge start operator $\star$
preserves $\Lambda^4(V)$ which splits therefore as $\Lambda^4(V)=\Lambda^4_{+}(V) \oplus \Lambda^4_{-}(V)$ into the $\pm$-eigenspaces of $\star$. Our discussion is divided into several
steps.
\subsection{Self-dual $4$-forms}

Let us pick $T$ in $\Lambda^4_{+}(V)$. Recall that in this case $Z_T=Z^{+}_T$ and consider the symmetric and
traceless operator $\mu_T : \SS^{+} \to \SS^{+}$. Let $\sigma_T=\{\lambda_q, 1 \le q \le p\}$ be the non-zero
part of the spectrum of $\mu_T$ where we assume the eigenvalues $\lambda_q, 1 \le q \le p$ to be pairwise
distinct and where we denote their multiplicities by $m_q, 1 \le q \le p$. Therefore we obtain a splitting
\begin{equation} \label{startspin}
\SS^{+}=Z_T \oplus \SS_1 \oplus \ldots \SS_p
\end{equation}
where $\SS_q$ are the eigenspaces of $\mu_T$ corresponding to the eigenvalues $\lambda_q, 1 \le q \le p$.  Our aim here is to examine the splitting of $\Lambda^2(V)$ induced by \eqref{startspin} and to relate
it directly to the form $T$. We need now to recall the following simple fact, which essentially exploits the squaring isomorphism in $8$-dimensions.
\begin{lema} \label{sqskew}
Let $x,y$ belong to $\SS^{+}$, and let $x \wedge y$ in $\Cl_8$ be given as
$$(x \wedge y) \psi = \langle \psi, x \rangle y - \langle \psi, y \rangle x $$
for all $\psi$ in $\SS$. Then : \\
(i)   $x \wedge y$ belongs to $\Cl_8^0 \cap \Cl_8^{+}$ and $(x \wedge y)^t=-x \wedge y$.\\
(ii)  if moreover $Tx=\lambda_1 x$ and $Ty=\lambda_2 y$ where $T$ belongs to $\Lambda^4_{+}$ then
        \begin{equation*}
        \begin{split}
        T(x \wedge y)T = \lambda_1 \lambda_2 x \wedge y, \quad \text{and} \quad
        T(x \wedge y)+(x \wedge y)T = (\lambda_1 +\lambda_2)x \wedge y.
        \end{split}
        \end{equation*}
(iii) under the assumptions in (ii), if $\lambda_1 = \lambda_2$ then $T(x \wedge y) = (x \wedge y)T = \lambda_1 x \wedge y $. \\
(iv)  if $x^{\prime}, y^{\prime}$ is another pair of spinors in $\SS^{+}$ then
        \begin{equation*}
        [x \wedge y, x^{\prime} \wedge y^{\prime}] = \langle x, y^{\prime} \rangle x^{\prime} \wedge y - \langle y, y^{\prime} \rangle x^{\prime} \wedge x
                                                   - \langle x, x^{\prime} \rangle y^{\prime} \wedge y + \langle x^{\prime}, y \rangle y^{\prime} \wedge x.
        \end{equation*}
(v)
        \begin{equation*}
        8 \langle x \wedge y, x^{\prime} \wedge y^{\prime} \rangle = \langle y, y^{\prime} \rangle \langle x, x^{\prime} \rangle -
                                                                     \langle x, y^{\prime} \rangle \langle y, x^{\prime} \rangle.
        \end{equation*}
\end{lema}

\begin{proof}
(i) is standard, see \cite{laws}. \\
We prove (ii) and (iii) at the same time.  For any $\psi$ in $\SS$ we have
\begin{equation*}
(x \wedge y) T \psi =  \langle T \psi, x \rangle y - \langle T \psi, y \rangle x
                    =  \langle T \psi, x \rangle y - \langle T \psi, y \rangle x
                    =  \lambda_1 \langle \psi, x \rangle y - \lambda_2 \langle \psi, y \rangle x
\end{equation*}
as $\langle T\psi, x \rangle = \langle \psi, Tx \rangle = \lambda_1 \langle \psi,x \rangle$ and similarly $\langle T \psi, y \rangle = \lambda_2 \langle \psi, y \rangle$. Moreover,
\begin{equation*}
T(x \wedge y)\psi  = \langle \psi, x \rangle T y - \langle \psi, y \rangle Tx
                   = \lambda_2 \langle \psi, x \rangle y - \lambda_1 \langle \psi, y \rangle x.
\end{equation*}
All claims in (ii) and (iii) follow now easily. The proof of (iv) is a straightforward direct computation involving only the definition of the exterior product of spinors. \\
(v) A direct computation based on the definition of the wedge product of spinors shows that the trace of the Clifford multiplication with the symmetric (in the sense of Lemma
\ref{trace}) element $(x \wedge y )(x^{\prime} \wedge y^{\prime})+(x^{\prime} \wedge y^{\prime})(x \wedge y)$ of $\Cl_8$ is given by
\begin{equation*}
-4 \Big[ \langle y, y^{\prime} \rangle \langle x, x^{\prime} \rangle - \langle x, y^{\prime} \rangle \langle y, x^{\prime} \rangle \Big].
\end{equation*}
The claim follows now by using Lemma \ref{trace}.
\end{proof}

For any $1 \le k,i,j \le p$ let us now define the spaces
\begin{equation*}
\begin{split}
E_k=&\{\alpha \in \Lambda^2 : T\alpha T=0, \ T \alpha+\alpha T=\frac{\lambda_k}{2}(1+\nu)\alpha \}\\
F_{ij}=&\{\alpha \in \Lambda^2 : T\alpha T=\frac{1}{2}\lambda_i \lambda_j (1+\nu)\alpha, \ T \alpha+\alpha T=\frac{\lambda_i+\lambda_j}{2}(1+\nu)\alpha \}\\
\iota_T^0=&\{ \alpha \in \Lambda^2 : \alpha T=T \alpha=0 \}.
\end{split}
\end{equation*}
Obviously  we have $F_{ij}=F_{ji}$. For notational convenience, we set $E=\bigoplus \limits_{k=1}^{p}E_k, F=\bigoplus \limits_{1 \le i \le j \le p}^{} F_{ij}$. Another related object is
\begin{defi}
The isotropy algebra $\iota_T$ of $T$ in $\Lambda^4$ is the subalgebra of $\mathfrak{so}(V)$ given by
\begin{equation*}
\{\alpha \in \mathfrak{so}(V) : [\alpha, T]=0 \}.
\end{equation*}
Here the Lie bracket is considered within the Lie algebra $\Cl_8$.
\end{defi}
\begin{pro} \label{split2}
The following hold :
\begin{itemize}
    \item[(i)]      there is an orthogonal, direct sum decomposition
                    \begin{equation*}
                        \Lambda^2=\iota_T^0 \oplus E \oplus F
                    \end{equation*}
    \item[(ii)]     we have the following string of isomorphisms:
                    \begin{equation*}
                    \begin{split}
                        &\iota^0_T  \cong \Lambda^2(Z_T) \\
                        &E_k  \cong Z_T \otimes \SS_k \\
                        &F_{kk}  \cong \Lambda^2(\SS_k) \\
                        &F_{ij}  \cong \SS_i \otimes \SS_j, i \neq j
                    \end{split}
                    \end{equation*}
    \item[(iii)]
                    \begin{equation*}
                        \iota_T=\iota_T^0 \oplus \bigoplus \limits_{k=1}^{p}F_{kk}
                    \end{equation*}
                    \end{itemize}
\end{pro}

\begin{proof}
We prove (i) and (ii) together. For any two vector sub-spaces $V,W$ of $\SS^{+}$ we denote by $V \hat{\otimes} W$ the inclusion of $V \otimes W$ into $\Lambda^2(\SS^{+})$. Letting now
$\iota_T^0, E_k, F_{ij}$ be the images of $\Lambda^2(Z_T), Z_T \hat{\otimes} \SS_k, \SS_i \hat{\otimes} \SS_j$ under the inverse of the linear isomorphism $\mu_{\vert (1+\nu)\Lambda^2}: (1+\nu)
\Lambda^2 \to \Lambda^2(\SS^{+})$ proves our claims by making use of Lemma \ref{sqskew}. \\
(iii) Pick $\alpha$ in $\iota_T$. Then $T\alpha=\alpha T$ hence $T^2\alpha+\alpha T^2=2T \alpha T$. It is easy to see that the operator $\{T^2, \cdot \}-2T\cdot T$ equals $0$ on
$\iota_T^0, \lambda^2_k 1_{E_k}$ on $E_k$ and $(\lambda_i-\lambda_j)^2 1_{E_{ij}}$ on $E_{ij}$ thus $\iota_T \subseteq \iota_T^0 \oplus \bigoplus \limits_{k=1}^{p}F_{kk}$. The reverse inclusion and therefore
the equality follows from the construction of the spaces $F_{kk}, 1 \le k \le p$ and Lemma \ref{sqskew}, (iii).
\end{proof}
In the Proposition above the fact that $V$ is $8$-dimensional,  which implies that $A^0 \cap \Cl_8^{+}=(1+\nu) \Lambda^2$ had been used in a crucial way. It is obviously valid on $\Lambda_{-}^4$ as well.
The block structure of the isotropy algebra of a form $T$ in $\Lambda^4_{+}$ has been already obtained in \cite{DHM}, by a slightly different method and under the additional assumption
that $T$ is a calibration on $V$. In this case, the work in \cite{DHM} gives a complete geometric description of the resulting orbits.
In order to understand the structure of the holonomy algebra of $T$ we need to have a look at the Lie algebraic features of the splitting above.
\begin{coro} \label{flie1}
Let $T$ in $\Lambda^4_{+}(V)$ with $Z_T \neq (0)$ be given. Then $F$ is a Lie sub-algebra of $\Lambda^2$ isomorphic with $\mathfrak{so}(Z_T^{\perp})$.
\end{coro}
\begin{proof}
From the construction of $F$ the Clifford multiplication map gives an isometry $\mu: (1 + \nu)F \to \Lambda^2(Z_T^{\perp})$ by Lemma \ref{sqskew}, (v). Moreover, this is a Lie algebra isomorphism by
(iv) of the same Lemma.
\end{proof}

\begin{lema} \label{Lie1}
Suppose that $Z_T \neq (0)$. We have :
\begin{itemize}
    \item[(i)]   $[\iota_T^0, \iota_T^0] = \iota_T^0, \ [\iota_T^0, E_k] = E_k, \ [\iota_T^0, F_{ij}] = 0$ for any $1 \le k,i,j \le p$.
    \item[(ii)]  $[E_i, E_j] = F_{ij}$ for $i \neq j$ and $[E_i, E_i] = F_{ii} \oplus \iota_T^0$
    \item[(iii)] $[E_i, F_{jk}] = 0$ if $i \neq j, k$ and  $[E_i, F_{ij}] = E_j$
    \item[(iv)]  $[F_{ij}, F_{kl}] = 0$ if $(i,j) \cap (k,l) = \emptyset $
    \item[(v)]   $[F_{ij}, F_{ik}] = F_{jk}$ if $i,j,k$ are mutually distinct
    \item[(vi)]  $[F_{ij},F_{ij}] = F_{ii} \oplus F_{jj}$ and $[F_{ii}, F_{ik}] = F_{ik}$ provided that $m_i \ge 2$.
\end{itemize}
\end{lema}

\begin{proof}
Follows directly from the general formula in (iii) of Lemma \ref{sqskew} after inspecting the various possibilities of combining factors in the splitting of $\SS^{+}$ as given by
\eqref{startspin}.
\end{proof}
We need now to establish the analogue of Lemma \ref{Lie1} when the form $T$ has $Z_T = (0)$. In this case $\iota_T^0 = E = (0)$ and using the same arguments as previously we get

\begin{lema} \label{Lie2}
Let $T$ be in $\Lambda^4_{+}$ with $Z_T = (0)$. We have :
\begin{itemize}
    \item[(i)]   $[F_{ij}, F_{kl}] = 0$ if $(i,j) \cap (k,l) = \emptyset $
    \item[(ii)]  $[F_{ij}, F_{ik}] = F_{jk}$ if $i,j,k$ are mutually distinct
    \item[(iii)] $[F_{ij},F_{ij}] = F_{ii} \oplus F_{jj}$ and $[F_{ii}, F_{ik}] = F_{ik}$ provided that $m_i \ge 2$.
\end{itemize}
\end{lema}
\subsection{Structure of the commutators}

We shall give in this section a simplified expression, relying on the particular dimension, for the
generating space of the even part of the fix algebra $\mathfrak{g}_T^{\star}$, where $T$ belongs to $\Lambda^4_{+}$.

\begin{lema} \label{com4}
Let $T$ belong to $\Cl_8^0 \cap \Cl_8^{+}$ satisfy $T^t=T$. We have
\begin{equation*}
-4[X \lrcorner T, Y \lrcorner T]=2T\alpha T+\frac{1}{4}L(T^2 \alpha+\alpha T^2)+4\vert T \vert^2 (1-\nu) \alpha
\end{equation*}
for all $X,Y$ in $V$, where $\alpha$ in $\Lambda^2(V)$ is given as $\alpha=X \wedge Y$.
\end{lema}

\begin{proof}
Let $a : \Lambda^2(V)  \to \Cl_8$ be defined by setting
$$ a(\alpha)=\sum \limits_{i=1}^{8}e_iT^2(e_i \lrcorner \alpha)$$
for some orthonormal frame $\{ e_i, 1 \le i \le 8 \}$. Lemma \ref{dcom} actually says that
\begin{equation*}
-4[X \lrcorner T, Y \lrcorner T]=2T\alpha T+a (\alpha)
\end{equation*}
and we need only work out a simpler expression for the operator $a$. We compute
\begin{equation*}
e_i T^2 (e_i \lrcorner \alpha)=\frac{1}{2}e_i(T^2 \alpha)e_i-\frac{1}{2}(e_iT^2e_i)\alpha
\end{equation*}
leading to $a(\alpha)=\frac{1}{2}L(T^2 \alpha)-\frac{1}{2}L(T^2)\alpha$. But $L(T^2\alpha)=\frac{1}{2}L(T^2 \alpha+\alpha T^2)$ since $[T^2, \alpha]$ is a $4$-form, and moreover
$LT^2=-8\vert T \vert^2 (1-\nu)$.
\end{proof}
\subsection{Computation of $\mathfrak{g}_T^{\star,0}$}

For a given $T$ in $\Lambda^4_{+}(\mathbb{R}^{8})$, we shall compute now the even part $\mathfrak{g}_T^{\star,0}$ of its holonomy algebra.
The main technical ingredient in this section is contained in the following observation.

\begin{lema} \label{power}
Let $T$ be in $\Lambda^4_{+}$. Then $\mathfrak{g}_{T^{2k+1}}^{\star} \subseteq \mathfrak{g}_T^{\star}$ for all $k$ in $\mathbb{N}$.
\end{lema}

\begin{proof} Let $\{e_i, 1 \le i \le 8\}$ be an orthonormal basis in $V$ and consider the partial Casimir operator $C_T : \Cl_8 \to \Cl_8$ given by
\begin{equation*}
    C_T=\sum \limits_{i=1}^{8} [e_i \lrcorner T, [e_i \lrcorner T, \cdot ]].
\end{equation*}
Obviously, $C_T$ preserves the algebra  $ \mathfrak{g}_T^{\star}$, that is $C_T(\mathfrak{g}_T^{\star}) \subseteq \mathfrak{g}_T^{\star}$. A straightforward computation actually shows that
\begin{equation*}
    C_T\phi=\sum \limits_{i=1}^{8} (e_i \lrcorner T)^2 \phi +\phi (e_i \lrcorner T)^2-2(e_i \lrcorner T)\phi (e_i \lrcorner T)
\end{equation*}
for all $\phi$ in $\Cl_8$. We shall now compute $C_T(X \lrcorner \phi)$ where $\phi$ belongs to $\Cl_8^0 \cap \Cl_8^{+}$ is such that $\phi^t=\phi$ (equivalently $\phi$ is in
$\Lambda^4_{+} \oplus \mathbb{R}(1+\nu)$). We compute
\begin{equation*}
    \begin{split}
    -8(e_i \lrcorner T) (X \lrcorner \phi)(e_i \lrcorner T)
    =&(e_iT-Te_i)(X\phi-\phi X)(e_i T -T e_i)\\
    =&(-Te_iX\phi-e_iT\phi X )(e_i T-Te_i)\\
    =&-e_i T\phi X e_i T+T e_iX \phi T e_i\\
    =&-e_i T\phi(-2 \langle e_i, X \rangle - e_i X)T + T(-2 \langle e_i, X \rangle - X e_i) \phi T e_i
    \end{split}
\end{equation*}
henceforth after summation we get
\begin{equation*}
    4 \sum \limits_{i=1}^{8} (e_i \lrcorner T) (X \lrcorner \phi)(e_i \lrcorner T)=-[X, T\phi T]-\sh L(T \phi) X T + \sh T X L(\phi T).
\end{equation*}
Now
\begin{equation*}
    \begin{split}
    4\sum \limits_{i=1}^{8}(e_i \lrcorner T)^2 =  \sum \limits_{i=1}^{8} (e_i T-T e_i )^2
    =  LT T-LT^2+8T^2+T LT
    =  8T^2 + 8 \vert T \vert^2 (1-\nu)
    \end{split}
\end{equation*}
as $LT=0$ and $LT^2=-8 \vert T \vert^2 (1-\nu) $. A short computation using the stability relations gives now
\begin{equation*}
T^2(X \lrcorner \phi)+(X \lrcorner \phi) T^2=-\sh(-T^2\phi X+X \phi T^2)
\end{equation*}
hence in the end we obtain
\begin{equation*} \label{interm}
\begin{split}
4 C_T(X \lrcorner \phi)= & -4(X \phi T^2-T^2 \phi X) + 2 [X,T\phi T]\\
                         & +L (T\phi)XT-TX L(\phi T) + 16 \vert T \vert^2 (X \lrcorner \phi)
\end{split}
\end{equation*}
for all $X$ in $V$ and where $\phi$ belongs to $\Lambda^4_{+} \oplus \mathbb{R}(1+\nu)$. In particular, for $\phi=T^k, k$ in $\mathbb{N}$ this yields
\begin{equation*}
\begin{split}
4C_T(X \lrcorner T^k)=&X \lrcorner T^{k+2}+32 \langle T^k, T \rangle X \lrcorner T +16 \vert T \vert^2 (X \lrcorner T^k)
\end{split}
\end{equation*}
for all $X$ in $V$,  where we have use that $L(T^{k+1})=-8 \langle T^k,T \rangle (1-\nu)$. By induction, given that $C_T$ preserves $\mathfrak{g}_T^{\star}$ and that the later contains
$\{ X \lrcorner T : X \in V \}$ we arrive at $\{X \lrcorner T^{2k+1} : X \in V \} \subseteq \mathfrak{g}_T^{\star}$, for all $k$ in $\mathbb{N}$ and our claim follows.
\end{proof}
\begin{lema} \label{comb0}
Let $\lambda_i$ where $1 \le i \le p$ belong to $\sigma_T$. If $dim_{\mathbb{R}}Z_T \neq 0$ and
\begin{equation*}
8\vert T \vert^{2k+1}=\lambda_i^{2(2k+1)}
\end{equation*}
holds for all $k$ in $\mathbb{N}$ then $dim_{\mathbb{R}}Z_T=6$ and
$\sigma_T=\{\lambda_i, -\lambda_i \}$ with multiplicities $(1,1)$, provided that $T \neq 0$.
\end{lema}

\begin{proof}
Since $16 \vert T \vert^{2k+1}=\sum \limits_{q=1}^p m_q \lambda_q^{2(2k+1)}$ by making use of Lemma \ref{trace}, the equation we have to solve becomes
\begin{equation} \label{mmain}
\sum \limits_{q=1}^p m_q \lambda_q^{2(2k+1)}=2\lambda_i^{2(2k+1)}
\end{equation}
for all $k$ in $\mathbb{N}$. We now divide by $\lambda_i^{2(2k+1)}$ and take the limit when $k \to \infty$. It follows that $\vert \lambda_q \vert \leq \vert \lambda_i|$ for all $1
\leq q \leq i$ and also that $\sum \limits_{\vert \lambda_q \vert=\vert \lambda_i \vert}^{}m_q=2$. It follows easily that $m_i=1$, otherwise we would have $m_i=2$ and further
$\sigma_T=\{\lambda_i\}$ by making use of \eqref{mmain}, which contradicts that $\mu_T$ is traceless. Therefore $-\lambda_i$ belongs to $\sigma_T$, with multiplicity $1$ and our
claim follows again from \eqref{mmain}.
\end{proof}

\begin{pro}\label{at2} Let $T$ belong to $\Lambda^4_{+}(V)$ with $dim_{\mathbb{R}}Z_T \neq 0,6$. Then
\begin{equation*}
(1+\nu) F  \oplus (1-\nu) \Lambda^2 \subseteq \mathfrak{g}_T^{\star, 0}.
\end{equation*}
\end{pro}

\begin{proof}
Making use of Lemma \ref{com4} we have that $\mathfrak{g}_T^{\star,0}$ contains the set
\begin{equation} \label{ggen}
\{2 T\alpha T+\frac{1}{4}L(T^2\alpha+\alpha T^2)+4\vert T \vert^2 (1-\nu)\alpha : \alpha \in \Lambda^2 \}.
\end{equation}
as this is just spanned by double commutators of elements in its generating set. From the above we find that $(1-\nu) \iota_T^0$ is contained in $\mathfrak{g}_T^{\star,0}$. Actually,
by using Lemma \ref{power} we have that $X \lrcorner T^{2k+1}$ belongs to $\mathfrak{g}_T^{\star}$ and therefore, after taking double commutators of such elements and using again Lemma
\ref{com4} we get that
\begin{equation*}
2T^{2k+1}\alpha T^{2k+1}+\frac{1}{4}L(T^{2(2k+1)}\alpha+\alpha T^{2(2k+1)})+4\vert T^{2k+1} \vert^2 (1-\nu) \alpha
\end{equation*}
belongs to $\mathfrak{g}_T^{\star,0}$ for any $\alpha$ in $\Lambda^2$. Now if $\alpha$ is in $E_i$, for some $1 \leq i \leq p$ we have $T\alpha T=0$ and an easy computation by
induction shows
$$ T^{2(2k+1)}\alpha+\alpha T^{2(2k+1)}=\frac{1}{2}\lambda_i^{2(2k+1)}(1+\nu)\alpha
$$
for all $k$ in $\mathbb{N}$. We are led eventually to having
$$(4\vert T^{2k+1}\vert^2-\frac{1}{2}\lambda_i^{2(2k+1)})(1-\nu)E_i$$
contained in $\mathfrak{g}_T^{\star}$ for all $1 \leq i \leq p$ and all $k$ in $\mathbb{N}$. \\
Now since $\sigma_T$ has not the form in Lemma \ref{comb0}, in other words $dim_{\mathbb{R}}Z_T \neq 6$, for each $1 \leq i \leq p$ the factor above will be non-vanishing for some $k$
in $\mathbb{N}$ whence $(1-\nu)E_i \subseteq \mathfrak{g}_T^{\star}$ whenever $1 \le i\le p$. Now taking commutators and using (ii) of Lemma \ref{Lie1} it follows that $(1 - \nu)
\Lambda^2 \subseteq \mathfrak{g}_T^{\star}$. But $L(T^2 \alpha+\alpha T^2)$ belongs to $(1-\nu)\Lambda^2$ for all $\alpha$ in $\Lambda^2$ hence we get from \eqref{ggen} that $\{T\alpha
T : \alpha \in \Lambda^2\}$ is contained in $\mathfrak{g}_{T}^{\star,0}$. Making use of the splitting in Proposition \ref{split2} this actually says that $(1+\nu)F \subseteq
\mathfrak{g}_{T}^{\star}$ and we have showed that
\begin{equation*}
(1+\nu)F \oplus (1-\nu) \Lambda^2 \subseteq \mathfrak{g}_T^{\star,0}.
\end{equation*}
\end{proof}

Therefore, when the set of spinors fixed by some self-dual $4$-form is not empty we can conclude, with one exception, that:

\begin{pro} \label{at0}
Let $T$ in $\Lambda^4_{+}(V)$ be given, and suppose that $dim_{\mathbb{R}}Z_{T} \neq 0, 6$. Then
    \begin{itemize}
    \item[(i)] $\mathfrak{g}_T^{\star,0}=(1+\nu) F  \oplus (1-\nu) \Lambda^2$
    \item[(ii)] $Z(\mathfrak{g}_T^{\star, 0})=(0)$.
    \end{itemize}
\end{pro}

\begin{proof}
(i) By Proposition \ref{at2},
it is enough to see that $\mathfrak{g}_T^{\star,0} \subseteq (1 + \nu) F  \oplus (1-\nu) \Lambda^2$ and this will be achieved by showing that $\mathfrak{g}_T^{\star,0}$ is orthogonal
to $(1+\nu)(\iota_0^{T} \oplus E)$. Indeed, by the definition of $Z_T$ we have $(X \lrcorner T) Z_T=0$ for all $X$ in $V$, therefore $\mathfrak{g}_T^{\star}=0$. We now pick $\phi$ in
$\mathfrak{g}_T^{\star,0}, x$ in $Z_T$ and $y$ in $\SS^{+}$. From the definition of $x \wedge y$ and $\phi x=0$ follows that
\begin{equation*}
\begin{split}
\phi (x \wedge y) \psi = &   \langle x, \psi \rangle \phi y \\
(x \wedge y) \phi \psi = & - \langle \phi \psi , y \rangle x = \langle \psi, \phi y \rangle x
\end{split}
\end{equation*}
for all $\psi$ in $\SS$. Since $\phi (x \wedge y)+\phi (x \wedge y)$ is a symmetric element in $\Cl_8$ after taking the trace we get
$$ Tr(\phi (x\wedge y) + (x \wedge y)\phi) = 2 \langle \phi y, x \rangle = -2 \langle y, \phi x \rangle = 0. $$
Now using Lemma \ref{trace} it follows that $\langle \phi (x \wedge y)+(x \wedge y) \phi, 1 \rangle = 0$ and since $\alpha (\phi^t) = -\phi$ we arrive at $\langle \phi, x \wedge y
\rangle = 0$. Because $\{ x \wedge y : x \in Z_T, y \in \SS$ spans $(1 + \nu)(\iota_T^0 \oplus E)$ it follows that $\mathfrak{g}_T^{\star ,0}$ is orthogonal to $(1+\nu)(\iota_T^0
\oplus E)$, hence contained in $(1+\nu)F$ and the claim follows.\\
(ii) follows eventually from Lemma \ref{Lie1}.
\end{proof}

We leave out now the case when $dim_{\mathbb{R}}Z_T=6$, to be treated further on, and look at
the situation when there are no non-zero fixed spinors which needs first a combinatorics Lemma.

\begin{lema} \label{combin}
Let $\lambda_i, \lambda_j$ belong to $\sigma_T$ with $1 \le i \neq j \leq p$. If $dim_{\mathbb{R}}Z_T=0$ and
$$ 4\vert T^{2k+1} \vert^2-\frac{1}{2}(\lambda_i^{2(2k+1)}+\lambda_j^{2(2k+1)})=(\lambda_i \lambda_j)^{2k}(4\vert T\vert^2-\frac{1}{2}(\lambda_i^2+\lambda^2_j)) $$
holds for all $k$ in $\mathbb{N}$ then $\sigma_T$ must be one of the following:
\begin{itemize}
    \item[(i)]   $\sigma_T=\{\lambda, -\lambda\}$
    \item[(ii)]  $\sigma_T=\{ \lambda, -\lambda, \mu, -\mu, \sqrt{\vert \lambda \mu \vert}, -\sqrt{\vert \lambda \mu \vert}\}$ with multiplicities $(1,1,1,1,2,2)$, for some
                 $\vert \lambda \vert \neq \vert \mu \vert$
    \item[(iii)] $\sigma_T=\{ \lambda, -\mu, \mu, \pm \sqrt{\vert \lambda \mu \vert}\}$ with multiplicities $(2,1,1,4)$.
\end{itemize}
\end{lema}

\begin{proof}
The proof is given in the Appendix for it is of rather technical nature and the reader may skip it at a first reading.
\end{proof}

\begin{pro} Let $T$ belong to $\Lambda^{4}_{+}$ such that $Z_T=(0)$. Then either
\begin{itemize}
    \item[(i)]
                \begin{equation*}
                \mathfrak{g}_T^{\star, 0}=A^0 \cong \mathfrak{so}(8) \oplus \mathfrak{so}(8)
                \end{equation*}
                or
    \item[(ii)]
                \begin{equation*}
                \mathfrak{g}_T^{\star,0}\cong \mathfrak{so}(8)
                \end{equation*}
                case which occurs when $T$ is a unipotent element of $\Cl_8^{+}$, in the sense that $T^2=\lambda (1+\nu)$ for some $\lambda >0$.
\end{itemize}
\end{pro}

\begin{proof}
Using Lemma \ref{power} we have that $[X \lrcorner T^{2k+1}, Y \lrcorner T^{2k+1}]$ belongs to $\mathfrak{g}_T^{\star,0}$ for all $X,Y$ in $V$ and any natural number $k$. Now Lemma
\ref{com4} implies that
\begin{equation*}
    2T^{2k+1}\alpha T^{2k+1}+\frac{1}{4}L(T^{2(2k+1)}\alpha+\alpha T^{2(2k+1)})+4\vert T^{2k+1} \vert^2 (1-\nu)\alpha
\end{equation*}
belongs to $\mathfrak{g}_T^{\star}$ for any $\alpha$ in $\Lambda^2(V)$ and any $k$ in $\mathbb{N}$. Let now $1 \le i , j \le p$ and pick $\alpha$ in $F_{ij}$. An easy computation by induction gives
\begin{equation*}
    \begin{split}
    & T^{2k+1}\alpha T^{2k+1}=\frac{1}{2}(\lambda_i \lambda_j)^{2k+1}(1+\nu)\alpha\\
    & T^{2(2k+1)}\alpha+\alpha T^{2(2k+1)}=\frac{1}{2}(\lambda_i^{2(2k+1)}+\lambda_j^{2(2k+1)})(1+\nu)\alpha
    \end{split}
\end{equation*}
hence
\begin{equation} \label{power2}
    (\lambda_i \lambda_j)^{2k+1} (1+\nu) \alpha+\biggl [ 4\vert T^{2k+1} \vert^2-\frac{1}{2}(\lambda_i^{2(2k+1)}+\lambda_j^{2(2k+1)})\biggr ](1-\nu)\alpha
\end{equation}
belongs to $\mathfrak{g}_T^{\star,0}$ for all $\alpha$ in $F_{ij}$ and all $k$ in $\mathbb{N}$. But in the same time, for $k=0$, we have
\begin{equation*}
    \lambda_i \lambda_j (1+\nu)\alpha+\biggl [ 4\vert T \vert^2-\frac{1}{2}(\lambda_i^{2}+\lambda_j^2)\biggr ](1-\nu)\alpha
\end{equation*}
belongs to $\mathfrak{g}_T^{\star, 0}$ for all $\alpha$ in $F_{ij}$ hence
\begin{equation*}
    \biggl [ 4\vert T^{2k+1} \vert^2-\frac{1}{2}(\lambda_i^{2(2k+1)}+\lambda_j^{2(2k+1)})-(\lambda_i \lambda_j)^{2k}(4\vert T\vert^2-\frac{1}{2}(\lambda_i^2+\lambda^2_j)) \biggr ] (1-\nu)\alpha
\end{equation*}
is in $\mathfrak{g}_T^{\star,0}$ whenever $\alpha$ belongs to $F_{ij}$ and $k$ is in $\mathbb{N}$. If $T$ is not an unipotent element, nor $\sigma_T$ has the form in (ii) or (iii) of Lemma \ref{combin}
we use Lemma \ref{combin} to see that the factor above will be non vanishing for
some $k$. This leads to $(1-\nu)F_{ij} \subseteq \mathfrak{g}_T^{\star, 0}$ and further from \eqref{power2} we get that $(1+\nu)F_{ij} \subseteq \mathfrak{g}_T^{\star, 0}$, since $0 \notin \sigma_T$.
Therefore when the spectrum of $\mu_T$ is not as in (ii) or (ii) of Lemma \ref{combin} our claim follows easily from the above, while the case when $T$ is unipotent is covered by Theorem \ref{uni}. \\
It remains now to treat the remaining two cases both of which to be worked out directly from \eqref{power2}.\\
(a) $\sigma_T=\{\lambda, -\lambda, \mu, -\mu, \sqrt{\vert \lambda \mu\vert},
-\sqrt{\vert \lambda \mu\vert}\}$ with multiplicities $(1,1,1,1,2,2)$.\\
In this case we have
\begin{equation*}
    4 \vert T \vert^{2k+1}=\frac{1}{2}(\lambda^{2(2k+1}+\mu^{2(2k+1)})+\vert \lambda \mu \vert^{2k+1}
\end{equation*}
for all $k$ in $\mathbb{N}$. We shall label the distinct eigenvalues as $\lambda_i, 1 \le i \le 6$ in the order they are listed in $\sigma_T$. For any $\alpha$ in $F_{i5}, 1 \le i \le 6$ we have then from
\eqref{power2}
\begin{equation*}
    2(\lambda_i \sqrt{\vert \lambda \mu \vert})^{2k+1}(1+\nu)\alpha+\biggl [ (\lambda^{2(2k+1)}+\mu^{2(2k+1)}) +\vert \lambda \mu \vert^{2(2k+1)}-\lambda_i^{2(2k+1)} \biggr ](1-\nu)\alpha
\end{equation*}
for all $k$ in $\mathbb{N}$. w.l.o.g we may also assume that $\vert \lambda \vert <\vert \mu \vert$. Now if $i=1,2$ or $i=5,6$ dividing by $\mu^{2(2k+1)}$ and taking the limit with
$k \to \infty$ yields $(1-\nu)F_{i5} \subseteq \mathfrak{g}_T^{\star,0}$ whence the same holds for $(1+\nu) F_{i5}$. Now if $i=3,4$ we divide by $(\mu \sqrt{\vert \lambda \mu \vert})^{2k+1}$ and make
$k \to \infty$ to arrive at $(1+\nu)F_{i5}$ contained in $\mathfrak{g}_T^{\star,0}$ and from there to $(1-\nu) F_{i5} \subseteq \mathfrak{g}_T^{0}$. Now taking commutators and using Lemma
\ref{Lie2} it follows that $(1\pm \nu)F_{ij}$ is contained in $\mathfrak{g}_T^{\star,0}$ hence the latter contains $A^0$ and our claim follows.\\
(b) $\sigma_T=\{ \lambda, -\mu, \mu, \pm \sqrt{\vert \lambda \mu \vert}\}$ with multiplicities $(2,1,1,4)$.\\
Here we label the (distinct) eigenvalues of $\mu_T$ by $\lambda_i, 1 \le i \le 5$, in their order of appearance and note that from the traceless of $\mu_T$ we must have
$\vert \lambda\vert=4 \vert \mu \vert > \vert \mu \vert $. Apart from this differences the proof of (a) continues to hold without any change. \\
We have exhausted all possibilities and therefore our claim is finally proved.
\end{proof}

We conclude this section with giving the full description of the holonomy algebras of forms $T$ in $\Lambda^4_{+}$ with $Z_T=(0)$.

\begin{teo} \label{halg1}
Let $T$ in $\Lambda^4_{+}$ be given and suppose moreover that $Z_T = (0)$. Then either
\begin{equation*}
\mathfrak{g}_T^{\star} = A \cong \mathfrak{so}(8,8)
\end{equation*}
or
\begin{equation*}
\mathfrak{g}_T^{\star} \cong \mathfrak{so}(8,1).
\end{equation*}
The latter case occurs when $T$ is an unipotent element of $\Cl_8^{+}$, that is $T^2=\lambda (1+\nu)$ for some $\lambda >0$.  In both cases the fix algebra is perfect, that is $\mathfrak{h}^{\star}_T=\mathfrak{g}_T^{\star}$.

\end{teo}

\begin{proof}
Because $[\mathfrak{g}_T^{\star,0}, \mathfrak{g}_T^{\star, 1}] \subseteq \mathfrak{g}_T^{\star,1}$, when $\mathfrak{g}_T^{\star,0}=A^0$ we obtain that $\mathfrak{g}_T^{\star,1}$ is an
invariant sub-space of $A^1$, w.r.t the adjoint representation of $A^0$ on $A^1$. Because this is irreducible (see Lemma \ref{propA}, (iii)) we find that $\mathfrak{g}_T^{\star,1}=A^1$
whence $\mathfrak{g}_T^{\star}=A$. The case of an unipotent element has been treated in Theorem \ref{uni}.
\end{proof}

\subsection{The full holonomy algebra when $Z_T \neq (0)$}

In order to have a complete description of holonomy algebras of self-dual $4$-forms in $8$-dimensions it remains to understand the odd part $\mathfrak{g}^{\star,1}_T$ of $\mathfrak{g}_T^{\star}$
for some $T$ in $\Lambda^4_{+}$ when $Z_T \neq (0)$. We recall that in this situation $\mathfrak{g}_T^{\star, 0 }$ has been computed in Proposition \ref{at0}. Let us now define
\begin{equation*}
Q=\{ \phi \in A^1 : T\phi+\phi T=0 \}.
\end{equation*}
We shall also work with the symmetric tensor product of spinors $(x,y) \to x \odot y$ where $x,y$ belong to $\SS^{+}$ which is defined in analogy with the exterior product
of spinors we saw before.

\begin{lema} \label{ocom}
Let $T$ belong to $\Lambda^{4}_{+}$. The following hold:
\begin{itemize}
    \item[(i)]  $Q = \{ \phi \in A^1 : T \phi = \phi T = 0 \}$.
    \item[(ii)] the map $(x,y) \to x \odot y$ extends to an isomorphism $\SS^{-} \otimes Z_T \to Q$.
\end{itemize}
\end{lema}

\begin{proof} (i) If $T\phi+ \phi T=0$ with $\phi$ in $A^1$, left multiplication with $\nu$ gives $T \phi-\phi T=0$, hence our claim, while using that $\nu T=T$ and
$\nu \phi+\phi \nu=0$. \\
(ii) It is easy to see from (i) that for any $\phi$ in $Q$ the map $\mu_{\phi}$ is a symmetric endomorphism of $\SS$ such that $\mu_{\phi}{\SS}^{-} \subseteq Z_T$. Details are very
similar to previous proofs and therefore left to the reader.
\end{proof}

\begin{pro}
Let $T$ belong to $\Lambda^4_{+}$ with $Z_T \neq (0)$. We have that $\mathfrak{g}_T^{\star ,1}=Q^{\perp}$.
\end{pro}

\begin{proof} A direct computation shows that
\begin{equation*}
\begin{split}
-[(1-\nu)\beta, X \lrcorner T]&=(\beta \wedge X)T+T(\beta \wedge X)+[X \lrcorner \beta, T]
\end{split}
\end{equation*}
for all $\beta$ in $\Lambda^2$ and all $X$ in $V$. Given that $(1-\nu)\Lambda^2 \subseteq \mathfrak{g}_T^{\star, 0}$ it follows easily that  $T \phi+\phi T$ belongs to
$\mathfrak{g}_T^{\star, 1}$ for all $\phi$ in $\Lambda^3$ and further that this actually holds for all $\phi $ in $A^1=\Lambda^3 \oplus \Lambda^7$. This is because
$T\Lambda^7+\Lambda^7 T$ just gives the generating set of
$\mathfrak{g}_T^{\star}$ since $\Lambda^7=\nu \Lambda^1$. \\
Therefore $\mathfrak{g}_T^{\star,1}$ contains the image of the symmetric operator $\{T, \cdot \} : A^1 \to A^1$ hence $Q^{\perp}$. Using Lemma \ref{ocom}, an argument completely similar to the one in the
proof of Proposition \ref{at0}, (i) leads to $Q^{\perp} \subseteq \mathfrak{g}_T^{\star,1}$ hence to the proof of the claim.
\end{proof}

Therefore our main result on holonomy algebras of self-dual $4$-forms with fixed spinors from this section is

\begin{teo}
Let $T$ be in $\Lambda^{4}_{+}$ with $dim_{\mathbb{R}}Z_T \neq 0,6$. Then the Clifford multiplication realises a Lie algebra isomorphism
$$\mu : \mathfrak{g}_T^{\star} \to \mathfrak{so}(8,8-dim_{\mathbb{R}}Z_T).$$
In particular, we must have $  \mathfrak{h}_T^{\star} = \mathfrak{g}_T^{\star} $.
\end{teo}

\begin{proof}
It is enough to see that  the map $(x,y) \to x \odot y$ gives an isomorphism ${\SS}^{-} \otimes Z_T^{\perp} \to Q^{\perp}$ and the rest follows by collecting the results above.
\end{proof}

It is also easy to see that under the assumptions above, $\mathfrak{g}_T^{\star}$ is a perfect Lie algebra.

\section{The case when $dim_{\mathbb{R}}Z_T=6$}
In this section we shall continue to work on an $8$-dimensional Euclidean vector space $(V^8, \langle \cdot, \cdot \rangle)$ which is furthermore supposed to be oriented, with
orientation form given by $\nu$ in $\Lambda^8(V)$. We will assume that $T$ in $\Lambda^4_{+}(V)$ satisfies $dim_{\mathbb{R}}Z_T=6$, and our primary aim will be to compute the algebra
$\mathfrak{g}_T^{\star}$.
As we have seen this situation cannot be covered only by the previous methods so we need more information about the structure of such forms. Let therefore $\sigma_T=\{\lambda_1,
\lambda_2\}$ be the non-zero part of the spectrum of $\mu_T : {\SS}^{+} \to {\SS}^{+}$ with multiplicities $(1,1)$ and let us also recall that $\lambda_1+\lambda_2=0$. We equally
recall that in this case the splitting of $\Lambda^2$ from Proposition \ref{split2} becomes
\begin{equation} \label{split6}
\Lambda^2=\iota_T^0 \oplus E_1 \oplus E_1 \oplus F_{12}
\end{equation}
and in particular $F$ is reduced to the $1$-dimensional component $F_{12}$. In what follows we shall use the normalisation $\lambda_1=1$ as it is clear that rescaling the generating
form leaves a holonomy algebra unchanged.
\subsection{Spinor $2$-planes} We start by recalling the following
\begin{defi} \label{acompx}
Let $(V^{8}, \langle \cdot, \cdot \rangle)$ be a Euclidean vector space. An almost Hermitian structure consists in a linear almost complex structure $J$ which is orthogonal w.r.t. the
scalar product $ \langle \cdot, \cdot \rangle$. If moreover $V$ is oriented, with orientation given by $\nu$ in $\Lambda^n(V)$, $J$ is positive if $\omega^4=\lambda \nu$ for some
$\lambda>0$ where $\omega=\langle J\cdot, \cdot \rangle $ is the so-called K\"ahler form of $(\langle \cdot, \cdot \rangle,J)$.
\end{defi}
In what follows we shall keep all previous notations and also recall the following well known fact, see \cite{dado} for instance.

\begin{pro} \label{2dimspc}
Let $L \subset \SS^{+}$ be any oriented $2$-dimensional sub-space of positive spinors. Then $L$ determines a unique positive almost Hermitian structure, say $J$, on $V$.
\end{pro}
For later use, and by sending again the reader to \cite{dado},
we mention that $J$ is constructed such that $(1+\nu)\omega=x_1 \wedge x_2$ for any oriented orthogonal basis $\{x_1,x_2\}$ in $L$ with the convention that $\vert x_1 \vert^2=\vert x_2
\vert^2=2$, where $\omega=\langle J\cdot, \cdot \rangle$. It is not difficult to see that the converse of Proposition \ref{2dimspc} also holds, in the sense that any compatible almost
Hermitian structure $J$ defines a $2$-dimensional sub-space $L$ of $\SS^{+}$ which is explicitly given as $L=Ker(\mu_{\omega}^2+16)$. For any compatible almost Hermitian structure $J$
we denote by $\lambda^4$ the underlying real bundle of the canonical line bundle of $J$. Explicitly, $\lambda^4=\{ \alpha \in \Lambda^4 : \alpha(J \cdot, J \cdot, \cdot,
\cdot)=-\alpha\}$ and if moreover $J$ is positive $\lambda^4$ is contained in $\Lambda_{+}^4$ (see \cite{redbook}). Note that if the contrary is not specified all forms are real valued
in this setting. We also recall that in presence of an almost Hermitian structure $\Lambda^2=\Lambda_0^2 \oplus \mathbb{R}\omega $, an orthogonal, direct, sum and that $\lambda^{1,1}$
is the space of $J$-invariant $2$-forms on $V$.
\begin{lema} \label{ctype}
Let $L \subset \SS^{+}$ be two dimensional and oriented and let $J$ be the complex structure determined by $L$. Then $T$ in $\Lambda^2(V)$ satisfies $T L=0$ iff $T$ belongs to
$\lambda^{1,1}_0(V)$.
\end{lema}
\begin{proof}
This is an easy exercise
taking into account that from the construction of $J$ it follows
\begin{equation*} \label{cpxspin}
\begin{split}
JY x_1&=-Y x_2 \\
JY x_2&=Y x_1
\end{split}
\end{equation*}
for all $Y$ in $V$, where $\{x_1, x_2 \}$ is an oriented orthonormal basis in $L$.
\end{proof}
This essentially leads to having $\iota^0_T=\lambda^{1,1}_0$ fact to be used later on and which encodes the well-known special isomorphism $\mathfrak{su}(4)\cong \mathfrak{so}(6)$
\cite{laws}. Moving within the same circle of arguments also shows that
\begin{pro} \label{linebd}
Given any $2$-dimensional sub-space $L \subset \SS^{+}$, the map $(x,y) \to x \odot y$  extends to an isomorphism $S^2_0(L) \to \lambda^4$.
\end{pro}

\begin{pro} \label{k4}
Any $4$-form $T$ in $\Lambda^4_{+}$ with $dim_{\mathbb{R}}Z_T=6$ determines uniquely an $SU(4)$-structure. That is, there exists a compatible and positive almost  Hermitian structure
$J$ on $V$ such that $T$ belongs to $\lambda^4$. The isotropy algebra of $T$ is isomorphic to $\mathfrak{su}(4)$.
\end{pro}

\begin{proof}
Let $L=Z^{\perp}_T$ be the orthogonal complement of $Z_T$ in $\SS^{+}$. Since this is $2$-dimensional we get a positive almost Hermitian structure $J$ on $V$. Now $\mu_T$ is completely
determined by its restriction to $L$ which gives an element in $S^2_0(L)$ and the fact that $T$ belongs to $\lambda^4$ follows from Proposition \ref{linebd}. The claim concerning the
isotropy algebra follows from Proposition \ref{split2}, (iii) by making use of the above mentioned special isomorphism $\mathfrak{su}(4) \cong \mathfrak{so}(6)$.
\end{proof}

This shows how to construct examples of self-dual $4$-forms $T$ such that $Z_T$ is of dimension $6$. Similarly, from the classification of self-dual $4$-forms on $\mathbb{R}^8$
obtained in \cite{DHM} one can easily give a geometric description of the cases when $Z_T$ has smaller dimension, but for considerations of time and space we shall not present those
here.

\subsection{The holonomy algebra}

As a convenient intermediary object, we shall make use of the Lie sub-algebra $\mathfrak{g}_T^{\star, 2}$ of $\mathfrak{g}_T^{\star, 0} \subseteq A^0$ generated by the sub-set
\begin{equation*}
\{[X \lrcorner T, Y \lrcorner T] : X,Y \in V\}
\end{equation*}
of $A^0$. We point out that {\it{a priori}} $\mathfrak{g}_T^{\star,0} \neq \mathfrak{g}_T^{\star,2}$.

\begin{lema}
We have
\begin{equation*}
\mathfrak{g}_T^{\star,2}=(3+\nu)F \oplus (1-\nu) \iota^0_T
\end{equation*}
\end{lema}

\begin{proof}
Follows by a straightforward computation based on the fact that $\mathfrak{g}_T^{\star,2}$ is generated by the set given in \eqref{ggen} and on the equations defining the spaces $E$
and $F$.
\end{proof}

For notational convenience let us set $Q^1_T=\{X \lrcorner T: X \in V \}$ and also $Q^2_T=\{ X \lrcorner (T \alpha_{12}): X \in V \} $. Given that $T\alpha_{12}+\alpha_{12}T=0$ is
easily seen that $T\alpha_{12}$ belongs to $\Lambda^4_{+}(V)$ and hence $Q^k_T, k=1,2$ are both contained in $\Lambda^3(V)$.

\begin{lema} \label{com6}
The following hold:
\begin{itemize}
\item[(i)] $[(1-\nu)\iota_T^0, Q^1_T]=Q^1_T$ \item[(ii)] $[(3+\nu)F_{12}, Q^1_{T}] \subseteq Q^1_T \oplus Q^2_T$ \item[(iii)] $[Q^1_T, Q^2_T]=\mathfrak{g}^{\star, 2}$
\end{itemize}
\end{lema}

\begin{proof}
(i) If $\alpha$ belongs to $\iota^0_T$ and $X$ is in $V$, an easy computation using essentially that $T\alpha=\alpha T=0$ and the self-duality of $T$ yields
\begin{equation*}
[(1-\nu)\alpha, X \lrcorner T] = -2[X \lrcorner \alpha, T].
\end{equation*}
(ii) Recall that $T\alpha_{12}+\alpha_{12}T=0$ and again using the self-duality of $T$ we obtain after a short computation
\begin{equation*}
[(3+\nu)\alpha_{12}, X \lrcorner T] = 3[X, T\alpha_{12}]-2[X \lrcorner \alpha_{12}, T]
\end{equation*}
for all $X$ in $V$.\\
(iii) Because we also have $T\alpha_{12}T = -\frac{1}{2}(1 + \nu) \alpha_{12}$ it follows that $T^2 \alpha_{12} = \alpha_{12} T^2 = \frac{1}{2} (1 + \nu) \alpha_{12}$. Therefore, by
using mainly the stability relations and that $T$ is self-dual, we arrive after computing at some length at
\begin{equation*}
\begin{split}
-4[X \lrcorner T, Y \lrcorner (T\alpha_{12})] = (TXYT)\alpha_{12}+\alpha_{12}(TYXT) + \frac{1}{2}(1-\nu)(X \alpha_{12}Y+Y \alpha_{12}X)
\end{split}
\end{equation*}
whenever $X,Y$ belong to $V$. Now using \eqref{split6} it is easily seen that $[T\alpha T, \alpha_{12}]=0$ for all $\alpha$ in $\Lambda^2(V)$ hence our commutator becomes
\begin{equation*}
\begin{split}
-4[X \lrcorner T, Y \lrcorner (T\alpha_{12})]= - \langle X,Y \rangle (1 + \nu) \alpha_{12} + \frac{1}{2}(1-\nu)(X \alpha_{12}Y+Y \alpha_{12}X)
\end{split}
\end{equation*}
for all $X,Y$ in $V$. On the other hand, given that $\alpha_{12}$ induces a compatible almost complex structure $J$ on $V$ such that $\alpha_{12} = \langle J \cdot, \cdot\rangle$ we
actually have
\begin{equation*}
\begin{split}
X \alpha_{12}Y + Y \alpha_{12} X = & 2(X \wedge (Y \lrcorner \alpha_{12}) + Y \wedge (X \lrcorner \alpha_{12}))-2 \langle X, Y \rangle \alpha_{12}\\
                                 = & 2(X \wedge JY + Y \wedge JX)_0 - \langle X, Y \rangle \alpha_{12}
\end{split}
\end{equation*}
where the subscript indicates orthogonal projection onto $\Lambda^2_0$. Henceforth, our commutator reads finally
\begin{equation*}
\begin{split}
-4[X \lrcorner T, Y \lrcorner (T\alpha_{12})] = -\frac{1}{2} \langle X, Y \rangle (3 + \nu) \alpha_{12} + (1 - \nu)(X \wedge JY + Y \wedge JX)_0
\end{split}
\end{equation*}
for all $X,Y$ in $V$. Obviously, $(X \wedge JY+Y \wedge JX)_0$ belongs to $\lambda^{1,1}_0=\iota_T^0$ hence $[Q^1_T, Q^2_T] \subseteq \mathfrak{g}_T^{\star,2}$ and the equality follows
at once when using the linear isomorphism $S^2_0 \to \lambda^{1,1}_0, S \to SJ$.
\end{proof}
\begin{teo} \label{last6} Let $T$ belong to $\Lambda_{+}^4(V)$ satisfy $dim_{\mathbb{R}}Z_T=6$. Then $\mathfrak{g}_T^{\star}$ is isomorphic to
$\mathfrak{so}(6,2)$ and moreover $ \mathfrak{h}_T^{\star}= \mathfrak{g}_T^{\star} $.
\end{teo}
\begin{proof}
It now easy to infer from the above that $\mathfrak{g}_T^{\star}=\mathfrak{g}_T^{\star,2} \oplus Q^1_T \oplus Q^2_T$, therefore the claim on $ \mathfrak{g}_T^{\star} $ follows. The
proof is completed when recalling that the Lie algebra $\mathfrak{so}(6,2)$ has trivial center.
\end{proof}
The proof of Theorem \ref{main8} is now complete.  We end this section by pointing out that in the case above the Clifford multiplication map $\mu:\mathfrak{g}_T^{\star} \to
Hom(Z_T^{\perp}, \SS^{-})$ is no longer surjective.

$\\$ {\bf{Acknowledgements}}: During the preparation of this paper both authors were partially supported by University of Auckland grants. They warmly acknowledge the Centro di Giorgi in Pisa and S.
Salamon for giving them the opportunity to attend the meeting ''Geometry and topology"  where parts of this work were completed.  During an early stage of this paper the research of
P.-A.N. has been supported through the VW-foundation at the HU-Berlin, while the work of N. B. was funded by the
DFG. Special thanks go to I.Agricola, S.Chiossi, A.R.Gover and Th.Friedrich for a number of useful discussions.

\appendix

\section{Proof of Lemma \ref{combin}}
Since after use of Lemma \ref{trace}
\begin{equation*} \label{powernorm}
16 \vert T^{2k+1} \vert^2=\sum \limits_{q=1}^{p}m_q \lambda_q^{2(2k+1)}
\end{equation*}
for all $k$ in $\mathbb{N}$ we obtain further
\begin{equation} \label{main}
\frac{1}{4} \sum \limits_{q=1}^{p}m_q \lambda_q^{2(2k+1)}=\frac{1}{2}(\lambda_i^{2(2k+1)}+\lambda_j^{2(2k+1)})+(\lambda_i \lambda_j)^{2k}(4\vert T\vert^2-\frac{1}{2}(\lambda_i^2+\lambda^2_j))
\end{equation}
whenever $k$ belongs to $\mathbb{N}$. \\
{\bf{Case I}}:  $\vert \lambda_i \vert \neq \vert \lambda_j \vert $\\
To fix ideas let us assume that $\vert \lambda_i \vert < \vert \lambda_j \vert$. It follows that
\begin{equation*}
\sum \limits_{q=1}^{p}m_q \lim_{k \to \infty} (\frac{\vert \lambda_q \vert}{\vert \lambda_j \vert})^{2(2k+1)}=2
\end{equation*}
hence $\vert \lambda_q \vert \leq \vert \lambda_j \vert $ for all $1 \le q \le p$ and moreover
\begin{equation} \label{mult1}
\sum \limits_{\vert \lambda_q \vert=\vert \lambda_j \vert}m_q=2.
\end{equation}
Therefore
\begin{equation} \label{main1}
\frac{1}{4} \sum \limits_{\vert \lambda_q \vert \neq
\vert \lambda_j \vert}^{p}m_q \lambda_q^{2(2k+1)}=\frac{1}{2}\lambda_i^{2(2k+1)}+(\lambda_i \lambda_j)^{2k}(4\vert T\vert^2-\frac{1}{2}(\lambda_i^2+\lambda^2_j))
\end{equation}
whenever $k$ belongs to $\mathbb{N}$. Further on, after dividing by $(\lambda_i \lambda_j)^{2k+1}$ and taking the limit when $k \to \infty$ we get
\begin{equation*}
\frac{1}{4} \sum \limits_{\vert \lambda_q \vert \neq
\vert \lambda_j \vert}^{p} m_q \lim_{k \to \infty} (\frac{\lambda_q^2}{\vert \lambda_i \vert \vert \lambda_j\vert})^{2k+1}=
\frac{1}{\vert \lambda_i \vert \vert \lambda_j\vert}(4\vert T\vert^2-\frac{1}{2}(\lambda_i^2+\lambda^2_j))
\end{equation*}
leading to $\lambda_q^2 \leq \vert \lambda_i \vert \vert \lambda_j\vert$ for all $\vert \lambda_q \vert \neq \vert \lambda_j \vert$ and
\begin{equation} \label{info2}
\sum \limits_{\stackrel{\lambda_q^2=\vert
\lambda_i \vert \vert \lambda_j \vert}{ \vert \lambda_q \vert \neq \vert \lambda_j \vert}}m_q=\frac{4}{\vert \lambda_i \vert \vert \lambda_j\vert}(4\vert T\vert^2-\frac{1}{2}(\lambda_i^2+\lambda^2_j)).
\end{equation}
Therefore, when actualising \eqref{main1} by \eqref{info2} we get
\begin{equation} \label{main2}
\sum \limits_{\stackrel{\lambda_q^2 \neq \vert \lambda_i \vert \vert \lambda_j \vert}{\vert \lambda_q \vert \neq \vert \lambda_j \vert}}^{p} m_q\lambda_q^{2(2k+1)}=2\lambda_i^{2(2k+1)}
\end{equation}
for all $k$ in $\mathbb{N}$. We now divide by $\lambda_i^{2(2k+1)}$ and take the limit when $k \to \infty $ to find that
\begin{equation*}
\sum \limits_{\stackrel{\vert \lambda_q \vert \neq \vert \lambda_j \vert}{\lambda_q^2 \neq \vert \lambda_i \lambda_j \vert}}^{p}m_q \vert \frac{\lambda_q}{\lambda_i} \vert^{2(2k+1)}=2
\end{equation*}
which implies that $\vert \lambda_q \vert \leq \vert \lambda_i \vert $ provided that $\vert \lambda_q \vert \neq \vert \lambda_j \vert$ and $\lambda_q^2 \neq \vert \lambda_i \lambda_j \vert$, and also that
\begin{equation} \label{info3}
\sum \limits_{\stackrel{\vert \lambda_q \vert \neq \vert \lambda_j \vert }{\stackrel{\lambda_q^2 \neq \vert \lambda_i \lambda_j \vert}{\vert \lambda_q \vert=\vert \lambda_i \vert }}}^{p}m_q=2.
\end{equation}
When actualising \eqref{main2} by the equation above it follows that the set
$$\{ \lambda_q \in \sigma_T: \vert \lambda_q\vert \neq \vert \lambda_j \vert, \vert \lambda_q\vert \neq \vert \lambda_j \vert, \lambda_q^2 \neq \vert \lambda_i \lambda_j \vert\}$$
is actually empty, in other words $\vert \lambda_q \vert $ belongs to $\{ \vert \lambda_i \vert, \vert \lambda_j \vert, \sqrt{\vert \lambda_i \lambda_j \vert} \}$ whenever $\lambda_q$ belongs to $\sigma_T$. Moreover,
from \eqref{mult1} and \eqref{info3} combined with $\vert \lambda_i \vert \neq \vert \lambda_j \vert $ we get that $m_j \leq 2$ together with $m_i \leq 2 $. Also, $\sigma_T$ always contains eigenvalues
$\lambda_q$ with $\lambda_q^2=\vert \lambda_i \lambda_j \vert$; otherwise, $\sigma_T$ would be included in $\{ \pm \lambda_i, \pm \lambda_j \}$ and moreover \eqref{main2} would imply that either $m_i=2$ or
$m_1=1$ and $\sigma_T$ contains $-\lambda_i$ with multiplicity $1$,  both of which cannot hold
on an $8$-dimensional space. Now, by a case by case discussion we shall consider all possibilities.\\
(i) $m_i=m_j=2$.\\
From \eqref{mult1} and \eqref{info3} it follows that $-\lambda_i, -\lambda_j$ are not present in $\sigma_T$ hence the presence of eigenvalues $\lambda_q$ with
$\lambda_q^2=\vert \lambda_i \vert \vert \lambda_j \vert $, of which we can have at most $2$, leads to the cases
$\sigma_T=\{ \lambda_i, \lambda_j, \pm \sqrt{\vert \lambda_i \vert \vert \lambda_j \vert }\}$ or
$\sigma_T=\{ \lambda_i, \lambda_j, \sqrt{\vert \lambda_i \vert \vert \lambda_j \vert }, -\sqrt{\vert \lambda_i \vert \vert \lambda_j \vert }\}$. In all cases, after counting
possible multiplicities for the new eigenvalues and using that $\mu_T$ is traceless we arrive at $\vert \lambda_i \vert=\vert \lambda_j \vert$ a contradiction. \\
(ii) $m_i=1, m_j=2$. \\
In this case $-\lambda_i, \lambda_i$ belong to $\sigma_T$ each of which with multiplicity $1$. Because we have eigenvalues $\lambda_q$ in $\sigma_T$ with
$\lambda_q^2=\vert \lambda_i \lambda_j \vert$ we have either $\sigma_T=\{ \lambda_j, \lambda_i, -\lambda_i, \pm \sqrt{\vert \lambda_i \lambda_j \vert}\}$
or $\sigma_T=\{\lambda_j, \lambda_i, -\lambda_i,  -\sqrt{\vert \lambda_i \lambda_j \vert},  \sqrt{\vert \lambda_i \lambda_j \vert} \}$. The latter case cannot be retained
because $\mu_T$ is trace free whilst in the first case  the last eigenvalue has multiplicity $4$, as the ambient space is $8$-dimensional. \\
(iii) $m_i=2,m_j=1$.\\
This is completely similar to the case (ii) since $\lambda_i$ and $\lambda_j$, if no ordering is assumed, play dual roles. \\
(iv) $m_i=m_j=1$.\\
Here $\pm \lambda_i, \pm \lambda_j$ belong to $\sigma_T$, all of them being simple. The only ways to complete $\sigma_T$ are $\sigma_T=\{ \lambda_i , -\lambda_i, \lambda_j, -\lambda_j,
\pm \sqrt{\vert \lambda_i \lambda_j \vert} \}$ or $\sigma_T=\{ \lambda_i , -\lambda_i, \lambda_j, -\lambda_j,
-\sqrt{\vert \lambda_i \lambda_j \vert},  \sqrt{\vert \lambda_i \lambda_j \vert} \} $ but the first is quickly discarded after using that $\mu_T$ is trace free. Concerning the second one, again by the vanishing of the trace
and $\vert \lambda_i \vert \neq \vert \lambda_j \vert$ we find that the multiplicities of the last two eigenvalues can only equal $2$.
$\\$
{\bf{Case II}}: $\vert \lambda_i \vert=\vert \lambda_j \vert$\\
In this situation it is easily seen that our equation becomes
$\sum \limits_{q \neq i,j}m_q \lambda_q^{2(2k+1)}=\lambda_i^{4k} \sum \limits_{q \neq i,j}m_q\lambda_q^2$ for all $k$ in $\mathbb{N}$. Dividing by $\lambda_i^{4k}$ and taking the limit when
$k \to \infty$ we obtain as above that $\vert \lambda_q \vert \leq \vert \lambda_i \vert $ for all $q \neq i,j$ and also that
$\sum \limits_{\stackrel{\vert \lambda_q \vert=\vert \lambda_i \vert}{q \neq i,j}} m_q\lambda_q^2=\sum \limits_{q \neq i,j}m_q\lambda_q^2$. But this implies that
$\vert \lambda_q \vert=\vert \lambda_i \vert$ for all $q \neq i,j$ which is an impossibility for it would imply $\lambda_q=\pm \lambda_i$ hence $q=i,j$ as $\lambda_i+\lambda_j=0$.
Therefore $\sigma_T=\{\lambda_i, \lambda_j\}$ and the claim follows.

\end{document}